\everymath{\displaystyle}
\documentclass[preprint]{elsarticle}

\usepackage[utf8]{inputenc}
\usepackage{amsmath}
\usepackage{amsfonts}
\usepackage{amssymb}
\usepackage[per-mode=symbol,inter-unit-product=\cdot]{siunitx}
\usepackage{physics}
\usepackage{mathtools}
\usepackage{subcaption}
\usepackage{cases}
\usepackage[hypertexnames=false,
            colorlinks=true,
            pdfborderstyle={/S/U/W 0},
            citecolor=blue]{hyperref}
\usepackage[capitalise,noabbrev]{cleveref}
\usepackage{todonotes}
\usepackage{multicol}
\usepackage{multirow}
\usepackage{booktabs}
\usepackage[a4paper, total={6.5in, 8.66in}]{geometry}
\usepackage{comment}

\DeclareSIUnit\mmhg{mmHg}
\newcommand{\sialpha}{\kilo\gram\per\square\metre\per\second}

\newcommand{\fluid}{_\mathrm{f}}
\newcommand{\solid}{_\mathrm{s}}
\newcommand{\ale}{_\mathrm{ALE}}

\renewcommand{\hat}{\widehat}

\newcommand{\exex}{EFS$_1$}
\newcommand{\exim}{E$_1$FS$_\infty$}
\newcommand{\exin}[1]{E$_1$FS$_#1$}

\makeatletter
\def\ps@pprintTitle{%
 \let\@oddhead\@empty
 \let\@evenhead\@empty
 \def\@oddfoot{}%
 \let\@evenfoot\@oddfoot}
\makeatother

\journal{Journal of Computational Physics}

\begin{document}

\begin{frontmatter}

\title{A stable loosely-coupled scheme for cardiac electro-fluid-structure interaction}


\author[mox]{Michele Bucelli\corref{cor}}
\cortext[cor]{Corresponding author. \url{michele.bucelli@polimi.it}}

\author[labs]{Martin Geraint Gabriel}
\author[unibg]{Giacomo Gigante}
\author[mox,epfl]{Alfio Quarteroni}
\author[labs]{Christian Vergara}

\affiliation[mox]{
    organization={MOX, Dipartimento di Matematica, Politecnico di Milano},
    addressline={P.zza Leonardo da Vinci 32},
    city={20133 Milano},
    state={Italy}}
\affiliation[labs]{
    organization={LABS, Dipartimento di Chimica, Materiali e Ingegneria Chimica ``Giulio Natta'', Politecnico di Milano},
    addressline={P.zza Leonardo da Vinci 32},
    city={20133 Milano},
    state={Italy}}
\affiliation[unibg]{
    organization={Dipartimento di Ingegneria Gestionale, dell'Informazione e della Produzione, Università degli Studi di Bergamo},
    addressline={Via Salvecchio 19},
    city={24129 Bergamo},
    state={Italy}}
\affiliation[epfl]{
    organization={Mathematics Institute, EPFL},
    addressline={Av. Piccard, CH-1015 Lausanne},
    state={Switzerland (Professor Emeritus)}
}

\begin{abstract}
We present a loosely coupled scheme for the numerical simulation of the cardiac electro-fluid-structure interaction problem, whose solution is typically computationally intensive due to the need to suitably treat the coupling of the different submodels. Our scheme relies on a segregated treatment of the subproblems, in particular on an explicit Robin-Neumann algorithm for the fluid-structure interaction, aiming at reducing the computational burden of numerical simulations. The results, both in an ideal and a realistic cardiac setting,  show that the proposed scheme is stable at the regimes typical of cardiac simulations. From a comparison with a scheme with implicit fluid-structure interaction, it emerges that, while conservation properties are not fully preserved, computational times significantly benefit from the explicit scheme. Overall, the explicit discretization represents a good trade-off between accuracy and cost, and is a valuable alternative to implicit schemes for fast large-scale simulations.
\end{abstract}


\end{frontmatter}

\section{Introduction}

Mathematical and numerical modeling of the cardiac function can provide meaningful insight into physiology, as well as assist in the development of personalized treatment \cite{gray2018patient,niederer2019computational,quarteroni2017integrated,quarteroni2019mathematical,verzicco2022electro}. Several computational models of the human heart function have been proposed, often focusing on specific features of its function: electrophysiology \cite{arevalo2016arrhythmia,bucelli2021multipatch,del2022fast,gillette2021framework,piersanti2021modeling,romero2010effects,trayanova2011whole,vergara2016coupled}, electromechanics \cite{augustin2016patient,augustin2021computationally,baillargeon2014living,fedele2022comprehensive,gerach2021electro,gerbi2018monolithic,gurev2011models,levrero2020sensitivity,pfaller2019importance,piersanti20213d,regazzoni2022cardiac,salvador2021electromechanical,strocchi2020simulating,usyk2002computational}, hemodynamics \cite{chnafa2014image,karabelas2022global,terahara2022computational,this2020augmented,zingaro2021hemodynamics,zingaro2022geometric} or fluid-structure interaction \cite{brenneisen2021sequential,bucelli2022partitioned,cheng2005fluid,khodaei2021personalized,nordsletten2011fluid,zhang2001analysis}.

Usually, the remaining features are neglected or surrogated by means of simplified models.
While this approach can provide meaningful results in physiological \cite{gerach2021electro,karabelas2022global,regazzoni2022cardiac,zingaro2022geometric} as well as pathological scenarios \cite{prakosa2018personalized,salvador2021electromechanical}, the heart function is characterized by the coordinated interplay of different physical processes, each affecting every other in multiple ways \cite{quarteroni2019mathematical}. Therefore, models featuring fully coupled and three-dimensional representations of electrophysiology, active and passive mechanics and fluid dynamics have the potential of providing a very accurate description of the physics of the heart \cite{gerbi2018numerical,hosoi2010multi,santiago2018fully,sugiura2022ut,viola2020fluid,viola2022fsei}. Models of this kind have been employed e.g. in computational studies on ventricular assist devices \cite{bakir2018multiphysics} and for in-silico clinical trials on digital cohorts of bundle branch block patients \cite{viola2022gpu}.
However, this comes at a high price in terms of model complexity and computational cost. For this reason, electro-fluid-structure models for the cardiac function are seldom considered in the literature. In particular, modeling the fluid-structure interaction (FSI) effects between the cardiac muscle and the blood dynamics is computationally challenging \cite{bucelli2022partitioned,cheng2005fluid,einstein2010fluid,feng2019analysis,hirschhorn2020fluid,nordsletten2011fluid,zhang2001analysis}. This is due in part to the anisotropy and non-linearity of the constitutive laws of muscular tissue, but also to the so-called added-mass effect \cite{causin2005added,forster2007artificial}: since fluid and structure have similar densities, numerical methods must be carefully designed to avoid time instability while keeping under control computational costs. These issues become even more pressing since, in the cardiac context, the FSI model is driven by the active muscular contraction, in turn triggered by electrical excitation, adding to its overall complexity and computational burden.

In this framework, explicit yet stable FSI schemes are very attractive \cite{bukavc2013fluid,burman2022fully,burman2009stabilization,burman2014explicit,fernandez2015generalized,gigante2021choice,gigante2021stability,guidoboni2009stable,li2016stable}. We focus in particular on loosely coupled Robin-Neumann (RN) type schemes, as discussed in \cite{gigante2021choice,gigante2021stability}. Here, RN indicates the kind of interface conditions that are alternatively enforced at the interface between the fluid and the structure. Analysis and numerical experiments for such schemes have been performed in both idealized and realistic vascular cases, considering a linear, isotropic and passive material for the structure \cite{gigante2021choice,gigante2021stability}, highlighting their stability for suitable ranges of the interface Robin parameter and of the time step.

The aim of this paper is to introduce, for the first time, a complete loosely coupled segregated scheme for the cardiac electro-fluid-structure interaction (EFSI) problem. To this aim, we combine a loosely coupled scheme used so far to couple electrophysiology and FSI \cite{bucelli2022mathematical} with the explicit RN scheme proposed in \cite{gigante2021choice,gigante2021stability} for FSI. We assess the effectiveness of the proposed scheme in a cardiac context, where the structure is characterized by a non-linear constitutive law, anisotropy and active contraction driven by electrophysiology. Moreover, we consider all the four cardiac phases (systolic and diastolic phases, together with the two isovolumic phases). We compare the complete loosely coupled scheme for EFSI against an explicit-implicit scheme where fluid and structure are strongly coupled in a monolithic fashion. We also consider a hybrid scheme, in which a small number of RN iterations is performed. We compare the schemes in terms of both accuracy and computational efficiency. Numerical results indicate that the complete loosely coupled scheme is stable in time and allows for significant computational savings. We also highlight a good accuracy of the explicit algorithm when compared to the explicit-implicit one, except for a mass loss during the isovolumic phases, so that overall the latter allows to strike a compromise between accuracy and computational efficiency. Finally, our conclusions are supported by a numerical experiment performed for a realistic ventricular model.

The rest of the paper is structured as follows. In \cref{sec:models} we introduce the cardiac EFSI problem. In \cref{sec:numerics} we detail the numerical methods used, with reference in particular to the complete loosely coupled scheme, and in \cref{sec:results} we present numerical results and a comparison between the schemes under consideration. Finally, in \cref{sec:conclusions} we draw some conclusive remarks.

\section{Mathematical modeling of cardiac\texorpdfstring{\\}{}electrophysiology-fluid-structure interaction}
\label{sec:models}

Let us denote by $t \in [0, T]$ the independent time variable. We consider a time-dependent domain $\Omega(t) \subset \mathbb{R}^3$ representing a human left ventricle \cite{bucelli2022partitioned}. The domain is split into the fluid part $\Omega\fluid(t)$, representing the volume occupied by the blood inside the chamber, and a solid part $\Omega\solid(t)$, corresponding to the volume occupied by the cardiac muscle, such that $\overline{\Omega(t)} = \overline{\Omega\fluid(t) \cup \Omega\solid(t)}$ and $\Omega\fluid(t) \cap \Omega\solid(t) = \varnothing$. The two domains share an interface $\Sigma(t) = \partial\Omega\fluid(t) \cap \partial\Omega\solid(t)$. We denote by $\mathbf n(t)$ the unit vector normal to $\partial\Omega(t)$, outgoing from $\Omega(t)$, and to $\Sigma(t)$, outgoing from $\Omega\fluid(t)$. To keep the notation light, we shall henceforth drop the explicit dependence on time of the sets defined so far, e.g. we will denote $\Omega(t)$, at the time-continous level, simply by $\Omega$.

\begin{figure}
    \centering
    \begin{subfigure}{0.64\textwidth}
        \includegraphics[width=0.49\textwidth]{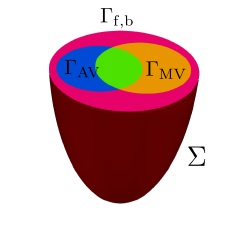}
        \includegraphics[width=0.49\textwidth]{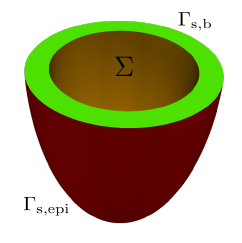}
        \caption{}
        \label{fig:domain}
    \end{subfigure}
    \begin{subfigure}{0.32\textwidth}
        \includegraphics[width=\textwidth]{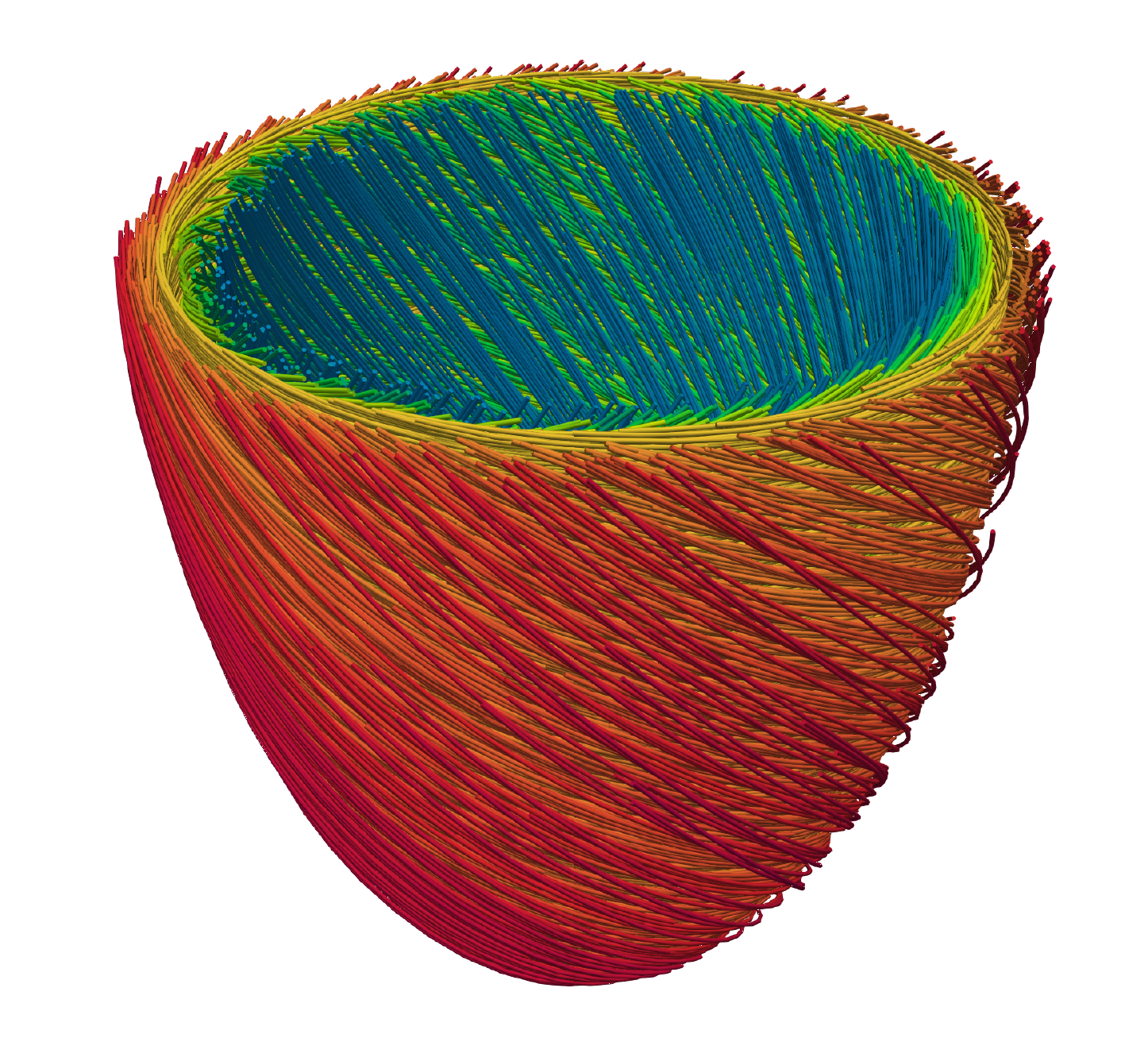}
        \caption{}
        \label{fig:fibers}
    \end{subfigure}
    \caption{(a) Computational domain $\Omega\fluid$ (left) and $\Omega\solid$ (right) of the idealized ventricle. Colors and labels denote the different portions of the boundary. (b) Streamline representation of the fiber field $\mathbf f_0$. Color is used to distinguish endocardium, myocardium and epicardium.}
\end{figure}

We denote by $\Gamma_\text{f,b}$ and $\Gamma_\text{s,b}$ the portion of boundary corresponding to the ventricular base on the fluid and structure domains, respectively. Moreover, we denote by $\Gamma_\text{s,epi}$ the epicardial surface (i.e. the outer wall), while $\Sigma$ corresponds to the endocardial surfaces on fluid and solid domains. There holds $\Gamma_\text{f,endo} = \Gamma_\text{s,endo} = \Sigma$. Finally, we denote by $\Gamma_\text{MV}$ and $\Gamma_\text{AV}$ two regions, possibly intersecting, of the fluid domain boundary representing the mitral and aortic valve orifices. The domain is represented in \cref{fig:domain}.

To track the motion of the domains, we introduce three fixed \textit{reference configurations} $\hat{\Omega}$, $\hat{\Omega}\fluid$, $\hat{\Omega}\solid$. We similarly denote with a hat the reference configuration for any of the previously defined boundary portions. The displacement of the time-dependent domains is expressed by the maps
\begin{align*}
    \mathcal{L}\solid &: \hat\Omega\solid \cross(0,T) \to \Omega\solid \qquad\qquad \Omega\solid(t) = \{\mathbf x = \mathcal{L}\solid(\hat{\mathbf x}, t)\;, \hat{\mathbf x}\in\hat\Omega\solid\}\;, \\
    \mathcal{L}\fluid &: \hat\Omega\fluid \cross(0,T) \to \Omega\fluid \qquad\qquad \Omega\fluid(t) = \{\mathbf x = \mathcal{L}\fluid(\hat{\mathbf x}, t)\;, \hat{\mathbf x}\in\hat\Omega\fluid\}\;.
\end{align*}
The precise definitions of $\mathcal{L}\solid$ and $\mathcal{L}\fluid$ are provided in the following sections.

The unknowns of our model are the following functions:
\begin{align*}
    & v:\hat{\Omega}\solid  \cross (0, T) \to \mathbb{R} & \text{transmembrane potential,} \\
    & \mathbf w:\hat{\Omega}\solid  \cross (0, T) \to \mathbb{R}^{N_\text{ion}^\mathbf{w}} & \text{gating variables,} \\
    & \mathbf z:\hat{\Omega}\solid  \cross (0, T) \to \mathbb{R}^{N_\text{ion}^\mathbf{z}} & \text{ionic concentrations,} \\
    & \mathbf s:\hat{\Omega}\solid  \cross (0, T) \to \mathbb{R}^{N_\text{act}} & \text{activation state variables,} \\
    & \mathbf d:\hat{\Omega}\solid  \cross (0, T) \to \mathbb{R}^3 & \text{solid displacement,} \\
    & \mathbf d\ale:\hat{\Omega}\fluid  \cross (0, T) \to \mathbb{R}^3 & \text{fluid domain displacement,} \\
    & \mathbf u:\Omega\fluid  \cross (0, T) \to \mathbb{R}^3 & \text{fluid velocity,} \\
    & p:\Omega\fluid  \cross (0, T) \to \mathbb{R} & \text{fluid pressure,}
\end{align*}
with $N_\text{ion}^\mathbf{w} = 12$, $N_\text{ion}^\mathbf{z} = 6$ and $N_\text{act} = 2$, according to the chosen ionic and force generation models (see \cref{sec:electrophysiology,sec:force-generation}).

\subsection{Fiber generation}
The cardiac tissue is characterized by the presence of fibers, that influence both its electrical and mechanical behavior \cite{eriksson2013influence,gil2019influence,piersanti2021modeling,roberts1979influence}. We account for their presence by defining at every point of $\hat{\Omega}\solid$ an orthonormal basis $\{\mathbf f_0, \mathbf s_0, \mathbf n_0\}$, representing the local directions of fibers, of fiber sheetlets and normal to fiber sheetlets, respectively. The basis is generated at every point by means of the algorithm presented in \cite{rossi2014thermodynamically}, as a preprocessing step. We refer the interested reader to \cite{piersanti2021modeling} for a detailed review of fiber generation methods for the whole heart. \Cref{fig:fibers} reports the generated fiber field on the idealized ventricle.

\subsection{Electrophysiology}
\label{sec:electrophysiology}

Electrophysiology models the evolution of the \textit{transmembrane potential}, i.e. the difference of potential $v$ between the intra- and extra-cellular spaces, as well as the evolution of ionic concentrations and ionic channels that determine the electrical excitation of cardiac cells \cite{collifranzone2014mathematical}. To this end, we use the monodomain equation, coupled with the ionic model by Ten Tusscher and Panfilov \cite{ten2006alternans}:
\begin{equation}
    \label{eq:electrophysiology}
    \begin{cases}
        \pdv{v}{t} - \div(\Sigma_\mathrm{m}\grad v) + I_\mathrm{ion}(v, \mathbf w, \mathbf z) = I_\mathrm{app}(\hat{\mathbf x}, t) & \text{in }\hat\Omega\solid \times (0, T)\;, \\
        \pdv{\mathbf w}{t} = \mathbf F_\mathrm{ion}^\mathbf{w}(v, \mathbf w) & \text{in } \hat\Omega\solid \times (0, T)\;, \\
        \pdv{\mathbf z}{t} = \mathbf F_\mathrm{ion}^\mathbf{z}(v, \mathbf w, \mathbf z) & \text{in } \hat\Omega\solid \times (0, T)\;, \\
        \Sigma_\mathrm{m}\grad v \cdot \mathbf n = 0 & \text{on }\partial\hat\Omega\solid \times (0, T)\;, \\
        v = v_0 & \text{in }\hat\Omega\solid \times \{0\}\;, \\
        \mathbf w = \mathbf w_0 & \text{in } \hat\Omega\solid \times \{0\}\;, \\
        \mathbf z = \mathbf z_0 & \text{in } \hat\Omega\solid \times \{0\}\;.
    \end{cases}
\end{equation}
In the above system, the first equation is the monodomain model, whereas the second and third equations express the ionic model in a compact form. We refer the interested reader to \cite{ten2006alternans} for the precise definition of $\mathbf F_\text{ion}^\mathbf{w}$, $\mathbf F_\text{ion}^\mathbf{z}$ and $I_\text{ion}$, as well as for the definitions of the components of $\mathbf w$ and $\mathbf z$. We point out that the vector $\mathbf z$ includes the intracellular concentration of calcium ions $[\text{Ca}^{2+}]_\text{i}$, that is relevant to the force generation model.

The tensor $\Sigma_\mathrm{m}$ expresses the anisotropic conduction properties of the myocardium. It is computed in terms of the fiber field as
\begin{equation*}
    \Sigma_\mathrm{m} =
    \sigma_\mathrm{m}^\mathrm{f}(\mathbf f_0 \otimes \mathbf f_0) +
    \sigma_\mathrm{m}^\mathrm{s}(\mathbf s_0 \otimes \mathbf s_0) +
    \sigma_\mathrm{m}^\mathrm{n}(\mathbf n_0 \otimes \mathbf n_0)\;,
\end{equation*}
where $\sigma_\mathrm{m}^\mathrm{f}$, $\sigma_\mathrm{m}^\mathrm{s}$ and $\sigma_\mathrm{m}^\mathrm{n}$ are conductivities in the fiber, sheetlet and normal directions, respectively \cite{piersanti2021modeling}.

Finally, the $I_\text{app}$ term in the monodomain equation models the ventricular stimulation by the Purkinje network in a simplified way, by applying a stimulus at three distinct locations on the endocardial wall \cite{piersanti2021modeling,regazzoni2022cardiac}.

We remark that we are neglecting the so-called geometry-mediated mechano-electric feedback effects \cite{salvador2022role}, that account for the fact that the electrical activation propagates in a moving domain. While relevant in pathological conditions, such effects have limited impact on simulations in sinus rhythm \cite{salvador2022role}. Nonetheless, the generalization of the proposed loosely coupled EFSI scheme to a model including mechano-electric feedback is straightforward.

Initial conditions $v_0$, $\mathbf w_0$ and $\mathbf z_0$ are obtained by running a single-cell electrophysiology simulation until a periodic limit cycle is reached, as detailed in \cite{regazzoni2022cardiac}.

\subsection{Force generation}
\label{sec:force-generation}

The state of contraction of cardiac cells is expressed at every point in $\hat\Omega\solid$ by the state vector $\mathbf s$. Its evolution is modeled by the ODE model proposed in \cite{regazzoni2018active}. Since the model features a very large number of variables, we use its reduced-order counterpart discussed in \cite{regazzoni2020machine}. The reduced model can be expressed as a system of ODEs defined at each point in $\hat\Omega\solid$:
\begin{equation}
    \label{eq:rdq18}
    \begin{cases}
        \pdv{\mathbf s}{t} = \mathbf{F}_\text{act}\left(\mathbf s, [\text{Ca}^{2+}]_\text{i}, \mathbf d\right) & \text{in } \hat\Omega\solid\times(0, T)\;, \\
        \mathbf s = \mathbf s_0 & \text{in } \hat\Omega\solid\times\{0\}\;.
    \end{cases}
\end{equation}
The generated active tension is then computed as a function of the contraction state as
\begin{equation*}
    T_\mathrm{act}(\mathbf s) = T_\mathrm{act,max} G(\mathbf s)\;
\end{equation*}
with $G(\mathbf s) \in [0, 1]$ and $T_\mathrm{act,max}$ the maximum generated contraction. We refer the interested reader to \cite{regazzoni2018active,regazzoni2020machine} for the precise definition of $\mathbf F_\text{act}$ and $G$. We remark that $\mathbf{F}_\text{act}$ depends on $\mathbf d$ accounting for the positive correlation between the stretch in the fiber direction and the generated force, coherently with the well-known Frank-Starling mechanism \cite{katz2010physiology,klabunde2011cardiovascular,opie2004heart}.

\subsection{Solid mechanics}

We define the map $\mathcal{L}\solid$ as
\begin{equation*}
    \mathcal{L}\solid(\hat{\mathbf x}, t) = \hat{\mathbf x} + \mathbf d(\hat{\mathbf x}, t)\;,
\end{equation*}
where $\mathbf d$ is the displacement field of the muscle, which is obtained as the solution of the elastodynamics equation in Lagrangian formulation \cite{ogden2013non}:
\begin{subnumcases}{\label{eq:elastodynamics}}
    \rho\solid \pdv[2]{\mathbf d}{t} - \div P\solid(\mathbf d, \mathbf s) = \mathbf 0 & in $ \hat\Omega\solid \times (0,T)\;,$ \\
    \mathbf d = \mathbf 0 & on $\hat\Gamma_\mathrm{s,b} \times (0, T)\;,$ \\
    P\solid(\mathbf d, \mathbf s)\mathbf n = -(\mathbf n \otimes \mathbf n)\left(K_\perp^\mathrm{epi}\mathbf d + C_\perp^\mathrm{epi}\pdv{\mathbf d}{t}\right) - (I - \mathbf n \otimes \mathbf n)\left(K_\parallel^\mathrm{epi}\mathbf d + C_\parallel^\mathrm{epi}\pdv{\mathbf d}{t}\right) & on $\hat\Gamma_\mathrm{s,epi} \times (0, T)\;,$ \label{eq:pericardium} \\
    \mathbf d = \mathbf d_0 & in $\hat\Omega\solid \times \{0\}\;,$ \\
    \pdv{\mathbf d}{t} = \mathbf 0 & in $\hat\Omega\solid \times \{0\}\;.$
\end{subnumcases}
In the above, $P\solid$ is the first Piola-Kirchhoff stress tensor, that accounts for both active and passive mechanical properties in the active stress framework \cite{ambrosi2012active}. It is decomposed additively as
\begin{equation*}
    P\solid(\mathbf d, \mathbf s) = P_\text{pas}(\mathbf d) + P_\text{act}(\mathbf d, \mathbf s)
\end{equation*}
into the passive contribution $P_\text{pas}(\mathbf d)$ and active contribution $P_\text{act}(\mathbf d, \mathbf s)$.

The passive stress tensor is defined in the hyperelastic framework as the derivative of a strain energy functional $\mathcal{W}$:
\begin{equation*}
    P_\text{pas}(\mathbf d) = \pdv{\mathcal{W}}{F}\;,
\end{equation*}
where $F = I + \grad\mathbf d$. We use the Guccione constitutive law for ventricular tissue \cite{guccione1991finite,regazzoni2022cardiac,usyk2002computational} with a penalization term for near-incompressibility. The constitutive law is non-linear and features anisotropy determined by the fiber field. See \ref{app:guccione} for more details.

The active contribution to the stress tensor is defined as \cite{regazzoni2022cardiac}
\begin{equation*}
    P_\mathrm{act}(\mathbf d, \mathbf s) = T_\mathrm{act}(\mathbf s)\frac{F\mathbf f_0\otimes\mathbf f_0}{\sqrt{I_\mathrm{4f}}}\;,
\end{equation*}
where $I_\mathrm{4f} = F\mathbf f_0 \cdot F\mathbf f_0$ measures the stretch along the fiber direction. We remark that $P_\mathrm{act}$ acts only in the direction of the fibers.

To find the initial displacement $\mathbf d_0$, we solve a quasi-static solid mechanics problem imposing a homogeneous endocardial pressure $p_0$ on $\hat\Gamma_\mathrm{s,epi}$ \cite{regazzoni2022cardiac}.

The condition \eqref{eq:pericardium} is a generalized visco-elastic Robin boundary condition that mimics the presence of the pericardial sac, a fluid-filled cavity that provides mechanical support, lubrication and protection from infections to the heart \cite{pfaller2019importance,regazzoni2022cardiac,strocchi2020simulating}.

Finally, we remark that the use of a homogeneous Dirichlet condition on the ventricular base is not consistent with physiology, and more sophisticated conditions should be applied \cite{piersanti20213d,regazzoni2022cardiac}. However, since our focus is on the numerical method for the FSI problem, we use a homogeneous Dirichlet condition for simplicity.

\subsection{Fluid domain displacement}

We model the motion of the fluid domain in the Arbitrary Lagrangian-Eulerian (ALE) framework \cite{donea1982arbitrary,duarte2004arbitrary,hughes1981lagrangian,nobile1999stability}. We introduce a fluid domain displacement field $\mathbf d\ale:\hat\Omega\fluid \to \mathbb{R}^3$ and define the mapping $\mathcal L\fluid$ as
\begin{equation*}
    \mathcal{L}\fluid(\hat{\mathbf x}, t) = \hat{\mathbf x} + \mathbf d\ale(\hat{\mathbf x}, t)\;.
\end{equation*}

The displacement $\mathbf d\ale$ is obtained by arbitrarily extending to $\hat\Omega\fluid$ the solid displacement on the interface, $\mathbf d|_{\hat\Sigma}$. We do so by means of a harmonic lifting operator:
\begin{equation}
    \label{eq:lifting}
    \begin{cases}
        -\Delta \mathbf d\ale = \mathbf 0 & \text{in } \hat\Omega\fluid \times (0, T)\;, \\
        \mathbf d\ale = \mathbf d & \text{on } \hat\Sigma \times (0, T)\;, \\
        \mathbf d\ale = \mathbf 0 & \text{on } \left(\hat\Gamma_\mathrm{f,b} \cup \hat\Gamma_\mathrm{AV} \cup \hat\Gamma_\mathrm{MV}\right) \times (0, T)\;.
    \end{cases}
\end{equation}

We define the ALE velocity $\mathbf u\ale$ as the time derivative of the ALE displacement, pushed forward to the current configuration:
\begin{equation*}
    \mathbf u\ale(\mathbf x, t) = {\pdv{\mathbf d\ale}{t}}(\mathcal{L}\fluid^{-1}(\mathbf x, t), t)\;.
\end{equation*}

\subsection{Fluid dynamics}
\label{sec:navierstokes}

We model the blood as an incompressible Newtonian fluid through Navier-Stokes equations in the ALE framework \cite{donea1982arbitrary,hughes1981lagrangian,quarteroni2017numerical}:
\begin{equation}
    \label{eq:navierstokes}
    \begin{cases}
        \rho\fluid\pdv{\mathbf u}{t} + \rho\fluid\left(\left(\mathbf u - \mathbf u\ale\right)\cdot\grad\right)\mathbf u - \div\sigma\fluid(\mathbf u, p) = \mathbf 0 & \text{in } \Omega\fluid \times (0, T)\;, \\
        \div\mathbf u = 0 & \text{in } \Omega\fluid \times (0, T)\;, \\
        \mathbf u = \mathbf 0 & \text{in } \Omega\fluid \times \{0\}\;, \\
        \mathbf u = \mathbf 0 & \text{on } \Gamma_\mathrm{f,b}\;,
    \end{cases}
\end{equation}
where $\rho\fluid$ is the fluid density,
\begin{gather*}
    \sigma\fluid(\mathbf u, p) = 2\mu\varepsilon(\mathbf u) - pI\;,\\
    \varepsilon(\mathbf u) = \frac{1}{2}\left(\grad\mathbf u + \grad\mathbf u^T\right)\;,
\end{gather*}
and $\mu$ is the dynamic viscosity of the blood. Suitable boundary conditions are imposed on $\Gamma_\text{MV}$ and $\Gamma_\text{AV}$ to model the opening and closing of the mitral and aortic valve, respectively, as detailed in \cref{sec:heartbeat-phases}.

\subsection{Fluid-structure interaction}

Besides the \textit{geometric coupling} expressed by \eqref{eq:lifting}, fluid and solid are coupled at the interface by imposing the continuity of velocity (\textit{kinematic coupling}) and of stresses (\textit{dynamic coupling}), expressing a no-slip condition and Newton's third law, respectively \cite{bazilevs2013computational}:
\begin{equation}
    \label{eq:interface}
    \begin{cases}
        \mathbf u = \pdv{\mathbf d}{t} & \text{on } \Sigma \times (0, T) \;,\\
        \sigma\fluid(\mathbf u, p)\mathbf n = \sigma\solid(\mathbf d, \mathbf s)\mathbf n & \text{on } \Sigma \times (0, T)\;,
    \end{cases}
\end{equation}
where $\sigma\solid(\mathbf d, \mathbf s)$ is the Cauchy stress tensor of the structure, related to $P(\mathbf d, \mathbf s)$ by
\begin{equation*}
    J \sigma\solid(\mathbf d, \mathbf s) = F P\solid(\mathbf d, \mathbf s)^T\;.
\end{equation*}
By taking a linear combination, the interface conditions \eqref{eq:interface} can be equivalently rewritten as follows \cite{badia2008fluid,gigante2021choice}:
\begin{equation}
    \label{eq:interface-rn}
    \begin{cases}
        \alpha\mathbf u + \sigma\fluid(\mathbf u, p)\mathbf n = \alpha\pdv{\mathbf d}{t} + \sigma\solid(\mathbf d, \mathbf s)\mathbf n & \text{on } \Sigma \times (0, T) \;,\\
        \sigma\fluid(\mathbf u, p)\mathbf n = \sigma\solid(\mathbf d, \mathbf s)\mathbf n & \text{on } \Sigma \times (0, T)\;,
    \end{cases}
\end{equation}
with $\alpha > 0$ a suitable Robin coefficient.

\subsection{Modeling of the four heartbeat phases}
\label{sec:heartbeat-phases}

Cardiac valves open and close passively to prevent reverse flow, and determine four distinct phases of the heartbeat \cite{katz2010physiology,klabunde2011cardiovascular,opie2004heart}. Focusing on the left heart, the phases are as follows:
\begin{enumerate}
    \item isovolumetric contraction: both the mitral and aortic valves are closed, and the ventricle starts to contract. This leads to a rapid increase in ventricular pressure, without any variation in ventricular volume. As soon as the ventricular pressure becomes larger than the pressure in the aorta, the aortic valve opens;
    \item ejection: blood is ejected from the ventricle into the aorta, leading to a decrease in ventricular volume. The mitral valve is closed, and the aortic valve is open. As soon as the flow through the aortic valve becomes null or negative, it closes;
    \item isovolumetric relaxation: both valves are again closed, and the ventricle starts relaxing. Ventricular pressure reduces, while ventricular volume stays constant. When the ventricular pressure becomes smaller than the atrial pressure, the mitral valve opens;
    \item filling phase: the mitral valve is open and the aortic valve is closed. Blood flows from the atrium into the ventricle, whose volume increases. When the flow through the mitral valve becomes null or negative, it closes.
\end{enumerate}
Isovolumetric contraction and ejection form the systolic phase, during which the ventricle contracts, whereas isovolumetric relaxation and filling form the diastolic phase. In order to model a full heartbeat, all of these phases must be captured adequately. We point out that isovolumetric phases pose significant modeling challenges for computational fluid dynamics simulations \cite{chnafa2014image,this2020augmented,zingaro2022modeling}, while FSI models such as the proposed one can deal with them naturally \cite{bucelli2022mathematical}.

We use switching boundary conditions on $\Gamma_\text{MV}$ and $\Gamma_\text{AV}$ to model opening and closing of the mitral and aortic valve, as done in \cite{bucelli2022partitioned}. No-slip conditions $\mathbf u = \mathbf 0$ are used for closed valves, the open mitral valve is modeled through a Neumann boundary condition, $\sigma(\mathbf u, p)\mathbf n = -p_\text{MV}\mathbf n$, and the open aortic valve is modeled through a resistance boundary condition,
\begin{equation*}
    \sigma\fluid(\mathbf u, p)\mathbf n = -\left(p_\text{AV}^0 + R_\text{AV} \int_{\Gamma_\text{AV}}\mathbf u\cdot\mathbf n\right)\mathbf n\;.
\end{equation*}
The choice of a resistance condition on the aortic valve allows to account for the typical evolution in time of the aortic pressure \cite{mark1998atlas}, at the same time preventing spurious reflections of pressure waves \cite{janela2010comparing,quarteroni2016geometric,vignon2010outflow}. Conversely, atrial pressure can be approximated as constant in time.

Valves are instantaneously switched from closed to open when the pressure upstream becomes larger than that downstream. Conversely, they are switched from open to close when the flowrate through them becomes null or reversed (i.e. when there is outflow through the mitral orifice or inflow through the aortic orifice). Thus, the opening and closure of valves is regulated by the numerical simulation and not prescribed a priori.


\section{Numerical discretization}
\label{sec:numerics}

Due to the large size of the problem, as well as its multiphysics and multiscale nature, the efficient numerical solution of the fully coupled electrophysiology-fluid-structure interaction (EFSI) system is a challenging task. One possible approach is based on a fully monolithic solver \cite{gerbi2018numerical}, where the coupled problem is discretized at each time step into a single large non-linear system. While very robust and stable, this approach requires the development of a dedicated solver and the use of suitable non-linear and linear solvers and preconditioners. Instead, one can choose to solve separately the different subproblems. This can be done maintaining a strong coupling between the problems, by means of subiterations \cite{santiago2018fully}, resulting in an implicit partitioned approach. However, iterative schemes of this kind might suffer from convergence issues, and can quickly become more computationally expensive than their monolithic counterpart \cite{bucelli2022partitioned}.

In this context, we propose a fully loosely coupled scheme in which all problems are solved only once per time step, and coupling terms are treated in an explicit fashion. This segregated approach, while possibly requiring the time step parameter to be sufficiently fine to guarantee stability, has the potential of being very computationally efficient. Moreover, the proposed method is inherently modular, and allows to flexibly choose spatial and temporal discretizations for the different equations.

While this approach is very common for electromechanics simulations \cite{augustin2021computationally,gerach2021electro,dede2020segregated,piersanti20213d,regazzoni2022cardiac}, the FSI coupling is usually treated in a strongly coupled way \cite{hirschhorn2020fluid,nordsletten2011fluid,santiago2018fully,watanabe2004multiphysics,zhang2001analysis}. Indeed, when the fluid and solid have comparable densities, as in biological applications, the so-called added-mass effect \cite{causin2005added} may lead to stability issues when a loosely coupled scheme is considered. Nonetheless, some loosely coupled FSI schemes have been investigated in the cardiovascular modeling literature \cite{burman2014explicit,burman2022fully,fernandez2015generalized,gigante2021choice,gigante2021stability}. In particular, in \cite{gigante2021choice,gigante2021stability} the authors investigated the stability properties of loosely coupled fluid structure interaction schemes based on Robin interface conditions, showing that for suitable choices of the Robin coefficients one may obtain a stable method. For the coupling of fluid and structure in our EFSI problem, we rely on a particular case of the Robin-based schemes, the loosely coupled Robin-Neumann (RN) scheme \cite{burman2022fully,fernandez2015generalized,gigante2021choice,gigante2021stability}.

We introduce a partition of the time domain $(0, T)$ into equally spaced sub-intervals $(t^{n}, t^{n+1})$, and denote by $\Delta t = t^{n+1} - t^n$ the width of each sub-interval. We denote with a superscript $n$ over any solution variable the time-discrete approximation of that solution variable at time $t^n$ (e.g. $\mathbf u^n \approx \mathbf u(t^n)$). We use finite differences for the time discretization of the subproblems. In the remainder of this section, we detail the proposed loosely coupled EFSI time discretization scheme (denoted with the abbreviation \exex{}), as well as an explicit-implicit (\exim{}) scheme in which the fluid-solid coupling is treated implicitly and monolithically \cite{bucelli2022partitioned,heil2004efficient,jodlbauer2019parallel,richter2015monolithic,wick2013solving,zhang2001analysis}, still maintaining explicit the coupling with the electrophysiology problem. We also present a scheme based on performing \num{2} fluid-structure subiterations (\exin{2}). The \exim{} scheme will be used as a reference for comparing numerical results. The considered schemes and the corresponding abbreviations are summarized in \cref{tab:scheme-abbreviations}.

\begin{table}
    \centering
    \begin{tabular}{c p{3.2cm} p{2cm} p{2.5cm}}
        \toprule
           \textbf{Abbreviation}
         & \textbf{Electromechanical coupling}
         & \textbf{Geometric coupling}
         & \textbf{FSI coupling} \\
         \midrule
         \exex{} & explicit & explicit & explicit RN \\
         \exim{} & explicit & explicit & implicit \\
         \exin{2} & explicit & explicit & \num{2} RN iterations \\
         \bottomrule
    \end{tabular}
    \caption{Summary of the abbreviations used to indicate the considered EFSI schemes. All schemes considered treat electrophysiology, force generation and fluid domain displacement in an explicit way.}
    \label{tab:scheme-abbreviations}
\end{table}

\subsection{Fully loosely coupled EFSI scheme (\texorpdfstring{\exex{}}{EFS1})}
\label{sec:explicit-scheme}

We detail in what follows the steps composing at each time step the fully loosely coupled FSI scheme (\exex{}). Given the solution up to time step $t^n$, in order to compute the solution at time $t^{n+1}$:

\begin{enumerate}
    \item Solve the electrophysiology model \eqref{eq:electrophysiology} with the following implicit-explicit (IMEX) scheme \cite{regazzoni2022cardiac}:
    \begin{enumerate}
        \item Solve the non-linear ionic model equations:
        \begin{equation}
            \label{eq:ionic-discrete}
            \begin{cases}
                \frac{\mathbf w^{n+1} - \mathbf w^n}{\Delta t} = \mathbf F_\mathrm{ion}^\mathbf{w}(v^n, \mathbf w^{n+1}) & \text{in } \hat\Omega\solid\;, \\[1em]
                \frac{\mathbf z^{n+1} - \mathbf z^n}{\Delta t} = \mathbf F_\mathrm{ion}^\mathbf{z}(v^n, \mathbf w^n, \mathbf z^n) & \text{in } \hat\Omega\solid\;;
            \end{cases}
        \end{equation}
        We remark that an implicit discretization is used for gating variables $\mathbf w$, whereas an explicit one is used for ionic concentrations $\mathbf z$;
        \item Solve the monodomain equation to compute $v^{n+1}$:
        \begin{equation}
            \label{eq:monodomain-discrete}
            \begin{cases}
                    \frac{v^{n+1} - v^n}{\Delta t} - \div(\Sigma_\mathrm{m}\grad v^{n+1}) + I_\mathrm{ion}(v^{n}, \mathbf w^{n+1}, \mathbf z^{n+1}) = I_\mathrm{app}^{n+1} & \text{in }\hat\Omega\solid\;, \\
                    \Sigma_\mathrm{m}\grad v^{n+1} \cdot \mathbf n = 0 & \text{on }\partial\hat\Omega\solid\;;
            \end{cases}
        \end{equation}
    \end{enumerate}
    \item Solve the time discretization of the force generation model \eqref{eq:rdq18}:
    \begin{equation}
        \label{eq:rdq18-discrete}
        \frac{\mathbf s^{n+1} - \mathbf s^n}{\Delta t} = \mathbf{F}_\text{act}(\mathbf s^n, [\text{Ca}^{2+}]_\text{i}^{n+1}, \mathbf d^n) \qquad \text{in } \hat\Omega\solid\;;
    \end{equation}
    \item Update the fluid domain solving \eqref{eq:lifting}:
    \begin{equation}
        \label{eq:lifting-discrete}
        \begin{cases}
            -\Delta \mathbf d\ale^{n+1} = \mathbf 0 & \text{in } \hat\Omega\fluid\;, \\
            \mathbf d\ale^{n+1} = \mathbf d^n & \text{on } \hat\Sigma\;, \\
            \mathbf d\ale^{n+1} = \mathbf 0 & \text{on } \hat\Gamma_\mathrm{f,b}\;;
        \end{cases}
    \end{equation}
    then set $\mathbf u\ale^{n+1} = \frac{\mathbf d\ale^{n+1} - \mathbf d\ale^n}{\Delta t}$ and compute the fluid domain at time $t^{n+1}$ as $\Omega\fluid^{n+1} = \mathcal{L}\fluid\left(\hat\Omega\fluid, t^{n+1}\right)$;
    \item Solve the time discretization of Navier-Stokes equations \eqref{eq:navierstokes} to compute $\mathbf u^{n+1}$ and $p^{n+1}$, with Robin boundary conditions on the fluid-solid interface:
    \begin{equation}
        \label{eq:navierstokes-discrete}
        \begin{cases}
            \begin{multlined}
            \rho\fluid\frac{\mathbf u^{n+1} - \mathbf u^{n}}{\Delta t} + \rho\fluid\left(\left(\mathbf u^{n} - \mathbf u\ale^{n+1}\right)\cdot\grad\right)\mathbf u^{n+1} -\div\sigma\fluid(\mathbf u^{n+1}, p^{n+1}) = \mathbf 0
            \end{multlined} & \text{in } \Omega\fluid^{n+1} \;, \\[1.5em]
            \div\mathbf u^{n+1} = 0 & \text{in } \Omega\fluid^{n+1} \;, \\
            \mathbf u^{n+1} = \mathbf 0 & \text{on } \Gamma_\mathrm{f,b}^{n+1}\;, \\
            \alpha\mathbf u^{n+1} + \sigma\fluid(\mathbf u^{n+1}, p^{n+1})\mathbf n^{n+1} = \alpha\frac{\mathbf d^{n} - \mathbf{d}^{n-1}}{\Delta t} + \sigma\solid(\mathbf d^n, \mathbf s^{n+1})\mathbf n^n & \text{on } \Sigma^{n+1}\;,
        \end{cases}
    \end{equation}
    endowed with suitable boundary conditions on $\Gamma_\text{MV}^{n+1}$ and $\Gamma_\text{AV}^{n+1}$ as described in \cref{sec:navierstokes}. We remark that interface conditions are computed using the solid displacement from previous time step, and that the advection term is treated in a semi-implicit way \cite{quarteroni2017numerical}, so that the resulting problem is linear;
    \item Solve the time discretization of the elastodynamics equation \eqref{eq:elastodynamics} to compute $\mathbf d^{n+1}$, with Neumann boundary conditions on the fluid-solid interface:
    \begin{equation}
        \label{eq:elastodynamics-discrete}
        \begin{cases}
            \rho\solid \frac{\mathbf d^{n+1} - 2\mathbf d^{n} + \mathbf d^{n-1}}{\Delta t^2} - \div P\solid(\mathbf d^{n+1}, \mathbf s^{n+1}) = \mathbf 0 & \text{in } \hat\Omega\solid\;, \\
            \mathbf d^{n+1} = \mathbf 0 & \text{on } \hat\Gamma_\mathrm{s,b}\;, \\
            \begin{multlined}P\solid(\mathbf d^{n+1}, \mathbf s^{n+1})\mathbf n \\[-1.00em] = -(\mathbf n \otimes \mathbf n)\left(K_\perp^\mathrm{epi}\mathbf d^{n+1} + C_\perp^\mathrm{epi}\frac{\mathbf d^{n+1} - \mathbf d^n}{\Delta t}\right) \\ - (I - \mathbf n \otimes \mathbf n)\left(K_\parallel^\mathrm{epi}\mathbf d^{n+1} + C_\parallel^\mathrm{epi}\frac{\mathbf d^{n+1} - \mathbf d^n}{\Delta t}\right)\end{multlined} & \text{on } \hat\Gamma\solid^\mathrm{epi}\;, \\[3.5em]
            \sigma\solid(\mathbf d^{n+1}, \mathbf s^{n+1})\mathbf n = \sigma\fluid(\mathbf u^{n+1}, p^{n+1})\mathbf n & \text{on } \Sigma\;.
        \end{cases}
    \end{equation}
    We point out that this problem is non-linear, due to the non-linearity of the constitutive law.
\end{enumerate}


\subsection{Explicit-implicit EFSI scheme (\texorpdfstring{\exim{}}{E1FS})}
\label{sec:explicit-implicit-scheme}

The explicit-implicit (\exim{}) scheme is based on treating explicitly the coupling of electrophysiology, force generation and mechanics, as well as the geometric FSI coupling, while treating the kinematic and dynamic FSI coupling in an implicit way. Given the solution up to time step $t^n$, in order to compute the solution at $t^{n+1}$, we repeat steps 1--3 as in the \exex{} scheme (\cref{sec:explicit-scheme}), replacing steps 4 and 5 with
\begin{enumerate}
    \setcounter{enumi}{3}
    \item Solve the time discretization of the FSI problem
    \begin{equation}
        \label{eq:fsi-implicit}
        \begin{cases}
            \rho\fluid\frac{\mathbf u^{n+1} - \mathbf u^{n}}{\Delta t} + \rho\fluid\left(\left(\mathbf u^{n} - \mathbf u\ale^{n+1}\right)\cdot\grad\right)\mathbf u^{n+1} -\div\sigma\fluid\left(\mathbf u^{n+1}, p^{n+1}\right) = \mathbf 0 & \text{in } \Omega\fluid^{n+1} \;, \\[1.5em]
            \div\mathbf u^{n+1} = 0 & \text{in } \Omega\fluid^{n+1} \;, \\
            \mathbf u^{n+1} = \mathbf 0 & \text{on } \Gamma_\mathrm{f,b}^{n+1}\;, \\
            \mathbf u^{n+1} = \frac{\mathbf d^{n+1} - \mathbf{d}^n}{\Delta t} & \text{on } \Sigma^{n+1}\;, \\
            \sigma\solid\left(\mathbf d^{n+1}, \mathbf s^{n+1}\right)\mathbf n = \sigma\fluid\left(\mathbf u^{n+1}, p^{n+1}\right)\mathbf n  & \text{on } \hat\Sigma\;, \\
            \rho\solid \frac{\mathbf d^{n+1} - 2\mathbf d^{n} + \mathbf d^{n-1}}{\Delta t^2} - \div P\solid(\mathbf d^{n+1}, \mathbf s^{n+1}) = \mathbf 0 & \text{in } \hat\Omega\solid\;, \\
            \mathbf d^{n+1} = \mathbf 0 & \text{on } \hat\Gamma_\mathrm{s,b}\;, \\
            \begin{multlined}P\solid\left(\mathbf d^{n+1}, \mathbf s^{n+1}\right)\mathbf n \\[-1.00em] = -(\mathbf n \otimes \mathbf n)\left(K_\perp^\mathrm{epi}\mathbf d^{n+1} + C_\perp^\mathrm{epi}\frac{\mathbf d^{n+1} - \mathbf d^n}{\Delta t}\right) \\ - (I - \mathbf n \otimes \mathbf n)\left(K_\parallel^\mathrm{epi}\mathbf d^{n+1} + C_\parallel^\mathrm{epi}\frac{\mathbf d^{n+1} - \mathbf d^n}{\Delta t}\right)\end{multlined} & \text{on } \hat\Gamma\solid^\mathrm{epi}\;,
        \end{cases}
    \end{equation}
    endowed with suitable boundary conditions on $\Gamma_\text{AV}^{n+1}$ and $\Gamma_\text{MV}^{n+1}$. This problem is non-linear due to the non-linearity of the solid constitutive law.
\end{enumerate}

We point out that both kinematic and dynamics FSI interface conditions are now treated implicitly. We will refer to this scheme as \exim{} for short.

\subsection{Explicit-hybrid EFSI scheme (\texorpdfstring{\exin{2}}{E1FS2})}

A hybrid approach between the \exex{} and the \exim{} schemes is obtained by introducing RN \cite{badia2008fluid,badia2009robin,gerardo2010analysis} subiterations with parameter $\alpha$ at each time step between fluid and structure problems, and performing \num{2} of such iterations. The steps 1--3 are the same as in the \exex{} scheme (\cref{sec:explicit-scheme}). Then, in place of the steps 4 and 5,
\begin{enumerate}
    \setcounter{enumi}{3}
    \item Setting $\mathbf u^{n+1}_{0} = \mathbf u^n$, $p^{n+1}_{0} = p^n$, $\mathbf d^{n+1}_{0} = \mathbf d^{n}$, iterate for $k = 0, 1$:
    \begin{enumerate}
        \item Solve the time discretization of the Navier-Stokes equations to compute $\mathbf u^{n+1}_{(k+1)}$ and $p^{n+1}_{(k+1)}$ by using structural displacement at previous iteration to prescribe Robin boundary conditions on the fluid-solid interface:
        \begin{equation*}
            \begin{cases}
                \rho\fluid\frac{\mathbf u^{n+1}_{(k+1)} - \mathbf u^{n}}{\Delta t} + \rho\fluid\left(\left(\mathbf u^{n} - \mathbf u\ale^{n+1}\right)\cdot\grad\right)\mathbf u^{n+1}_{(k+1)} -\div\sigma\fluid\left(\mathbf u^{n+1}_{(k+1)}, p^{n+1}_{(k+1)}\right) = \mathbf 0 & \text{in } \Omega\fluid^{n+1} \;, \\[1.5em]
                \div\mathbf u^{n+1}_{(k+1)} = 0 & \text{in } \Omega\fluid^{n+1} \;, \\
                \mathbf u^{n+1}_{(k+1)} = \mathbf 0 & \text{on } \Gamma_\mathrm{f,b}^{n+1}\;, \\
                    \alpha\mathbf u^{n+1}_{(k+1)} + \sigma\fluid\left(\mathbf u^{n+1}_{(k+1)}, p^{n+1}_{(k+1)}\right)\mathbf n^{n+1} = \alpha\frac{\mathbf d^{n+1}_{(k)} - \mathbf{d}^{n-1}}{\Delta t} + \sigma\solid\left(\mathbf d^{n+1}_{(k)}, \mathbf s^{n+1}\right)\mathbf n^n & \text{on } \Sigma^{n+1}\;;
            \end{cases}
        \end{equation*}
        \item Solve the time discretization of the elastodynamics equations to compute $\mathbf d^{n+1}_{(k+1)}$, using newly computed fluid velocity and pressure to provide Neumann conditions on the fluid-solid interface:
        \begin{equation*}
            \begin{cases}
                \rho\solid \frac{\mathbf d^{n+1}_{(k+1)} - 2\mathbf d^{n} + \mathbf d^{n-1}}{\Delta t^2} - \div P\solid\left(\mathbf d^{n+1}_{(k+1)}, \mathbf s^{n+1}\right) = \mathbf 0 & \text{in } \hat\Omega\solid\;, \\
                \mathbf d^{n+1}_{(k+1)} = \mathbf 0 & \text{on } \hat\Gamma_\mathrm{s,b}\;, \\[1em]
                \begin{multlined}P\solid\left(\mathbf d^{n+1}_{(k+1)}, \mathbf s^{n+1}\right)\mathbf n \\[-1.00em] = -(\mathbf n \otimes \mathbf n)\left(K_\perp^\mathrm{epi}\mathbf d^{n+1}_{(k+1)} + C_\perp^\mathrm{epi}\frac{\mathbf d^{n+1}_{(k+1)} - \mathbf d^n}{\Delta t}\right) \\ - (I - \mathbf n \otimes \mathbf n)\left(K_\parallel^\mathrm{epi}\mathbf d^{n+1}_{(k+1)} + C_\parallel^\mathrm{epi}\frac{\mathbf d^{n+1}_{(k+1)} - \mathbf d^n}{\Delta t}\right)\end{multlined} & \text{on } \hat\Gamma\solid^\mathrm{epi}\;, \\[4em]
                \sigma\solid\left(\mathbf d^{n+1}_{(k+1)}, \mathbf s^{n+1}\right)\mathbf n = \sigma\fluid\left(\mathbf u^{n+1}_{(k+1)}, p^{n+1}_{(k+1)}\right)\mathbf n  & \text{on } \hat\Sigma\;.
            \end{cases}
        \end{equation*}
    \end{enumerate}
    Then, set $\mathbf u^{n+1} = \mathbf u^{n+1}_{(2)}$, $p^{n+1} = p^{n+1}_{(2)}$ and $\mathbf d^{n+1} = \mathbf d^{n+1}_{(2)}$.
\end{enumerate}

We refer to this scheme as \textit{explicit-hybrid}, due to the fact that the iterative algorithm is arbitrarily truncated at \num{2} iterations, instead of checking for convergence, resulting in a hybrid approach between the explicit and implicit discretizations.

\subsection{Space discretization, non-linear and linear solvers}
The discretized-in-time problems introduced above are discretized using finite elements \cite{hughes2012finite,quarteroni2017numerical}. Independently of the scheme presented in the above sections, fluid and solid meshes are conforming at the interface $\Sigma$, and we stabilize the discretized Navier-Stokes equations using the SUPG-PSPG stabilization \cite{tezduyar2003stabilization}. Moreover, the ionic model \eqref{eq:ionic-discrete} and force generation model \eqref{eq:rdq18-discrete} are solved independently at each vertex of the computational mesh. For the ionic model, in particular, we adopt the \textit{ionic current interpolation} (ICI) approach \cite{krishnamoorthi2013numerical,pathmanathan2011significant}. The linear systems arising from the discretization of the monodomain equation \eqref{eq:monodomain-discrete} and of the fluid domain displacement problem \eqref{eq:lifting-discrete} are solved by means of the conjugate gradient (CG) method \cite{quarteroni2017numerical,saad2003iterative}, with an algebraic multigrid (AMG) preconditioner \cite{xu2017algebraic}.

Specifically to the \exex{} and \exin{2} schemes, the non-linear system arising from the solid mechanics discretization \eqref{eq:elastodynamics-discrete} is first linearized by means of Newton's method, and the resulting linear system is solved with GMRES \cite{saad2003iterative}, preconditioned using AMG. The block linear system arising from the discretization of Navier-Stokes equations \eqref{eq:navierstokes-discrete} is solved with GMRES with the SIMPLE preconditioner \cite{deparis2014parallel}, which in turn falls back onto AMG for the approximation of velocity and pressure diagonal blocks.

Instead, referring to the \exim{} scheme, for the solution of the FSI problem \eqref{eq:fsi-implicit}, we use a monolithic solver as presented in \cite{bucelli2022partitioned}, in which both fluid and solid equations are assembled in a single non-linear system. The latter is linearized with Newton's method and then solved with GMRES, using a block-lower triangular preconditioner that falls back onto SIMPLE and AMG for the fluid and structure submatrices.

\section{Numerical results}
\label{sec:results}

\begin{table}
    \centering
    \begin{tabular}{c c c c c c c c}
        \toprule
         & & \multicolumn{3}{c}{\textbf{Fluid}} & \multicolumn{3}{c}{\textbf{Structure}} \\
        \textbf{Mesh} & \textbf{Type} & \textbf{\# elem}. & \textbf{\# nodes} & \textbf{$h$ [\si{\milli\metre}]} & \textbf{\# elem}. & \textbf{\# nodes} & \textbf{$h$ [\si{\milli\metre}]} \\
        \midrule
        $\mathcal{M}_1$ & hex & \num{4684} & \num{5927} & \num{6.1} & \num{6612} & \num{8789} & \num{5.2} \\
        $\mathcal{M}_1^\text{EP}$ & hex & - & - & - & \num{52896} & \num{60459} & \num{2.6}\\
        $\mathcal{M}_2$ & hex &\num{13780} & \num{16669} & \num{4.0} & \num{22396} & \num{28117} & \num{3.3} \\
        $\mathcal{M}_3$ & hex &\num{32628} & \num{38429} & \num{3.0} & \num{51364} & \num{62589} & \num{2.5} \\
        \hline
        $\mathcal{M}_\text{R}$ & tet & \num{140644} & \num{157369} & \num{1.8} & \num{73860} & \num{89314} & \num{2.2} \\
        \bottomrule
    \end{tabular}
    \caption{Type of elements (hexahedra or tetrahedra), number of elements, number of nodes and average element diameter $h$ for the meshes considered on the prolate ellipsoid geometry ($\mathcal{M}_1$, $\mathcal{M}_1^\text{EP}$, $\mathcal{M}_2$ and $\mathcal{M}_3$) and for the realistic ventricle mesh ($\mathcal{M}_R$), for both the fluid and the structure domain.}
    \label{tab:mesh}
\end{table}

Numerical methods were implemented in \texttt{life\textsuperscript{x}} \cite{africa2022lifexcore,africa2022lifex,lifex}, a C++ high-performance computing library tailored at cardiac applications and based on the finite element core \texttt{deal.II} \cite{arndt2020dealii9.2,arndt2020dealii,dealii}.
In the following sections we report the results of numerical simulations using all the schemes presented in \cref{sec:numerics}, considering an idealized left ventricle described as a prolate ellipsoid (\cref{fig:domain}) and a realistic left ventricle model \cite{zygote}.
We report in \cref{tab:mesh} the discretization parameters of the the meshes under consideration. We compare the solutions obtained with the different schemes as well as their computational efficiency.

The values used for model parameters are reported in \ref{app:parameters}. Unless otherwise specified, simulations were run in parallel on 20 cores with CPUs Xeon E5-2640v4@2.4GHz, using the computational resources available at MOX, Mathematics Department, Politecnico di Milano.

\subsection{Solution indicators}
\label{sec:indicators}

One of the aims of the comparison among the schemes is to assess the loss of mass they feature, in particular during isovolumetric phases. To quantify this effect, we introduce two indices, the \textit{isovolumetric loss indices} (ILI), representing the relative variation of blood volume during isovolumetric phases:
\begin{equation*}
    \mathrm{ILI}_\mathrm{C} = \left|\frac{V_\mathrm{C,i} - V_\mathrm{C,f}}{\max\{V_\mathrm{C,i}, V_\mathrm{C,f}\}}\right|\qquad\qquad
    \mathrm{ILI}_\mathrm{R} = \left|\frac{V_\mathrm{R,i} - V_\mathrm{R,f}}{\max\{V_\mathrm{R,i}, V_\mathrm{R,f}\}}\right|\;,
\end{equation*}
wherein $V_\mathrm{C,i}$ and $V_\mathrm{C,f}$ are the volumes at the beginning and end of isovolumetric contraction, and $V_\mathrm{R,i}$ and $V_\mathrm{R,f}$ are the volumes at the beginning and end of isovolumetric relaxation. Optimal values for these two indices are $\mathrm{ILI}_\mathrm{C} = \mathrm{ILI}_\mathrm{R} = 0$, while positive values indicate that blood mass is not exactly preserved during the isovolumetric phases.

We also take into account the \textit{ejection fraction} EF and \textit{peak systolic pressure} $p_{\max}$, defined as
\begin{equation*}
    \mathrm{EF} = \frac{\text{EDV} - \text{ESV}}{\text{EDV}} \qquad\qquad
    p_{\max} = \max_{t \in (0, T)} \bar{p}(t)\;,
\end{equation*}
where $\bar{p}(t)$ is the ventricular average pressure at time $t$, and EDV and ESV are the end-diastolic and end-systolic volumes, i.e.
\begin{equation*}
    \text{EDV} = \max_{t \in (0, T)}V(t) \qquad\qquad
    \text{ESV} = \min_{t \in (0, T)}V(t)\;,
\end{equation*}
where $V(t)$ is the ventricular volume. Both $\mathrm{EF}$ and $p_{\max}$ have significant clinical relevance \cite{katz2010physiology,klabunde2011cardiovascular,opie2004heart}.

\subsection{Test A: on the stability of the loosely coupled scheme}
\label{sec:test-a}

\begin{figure}
    \centering
    \begin{subfigure}{0.4\textwidth}
        \includegraphics[width=\textwidth]{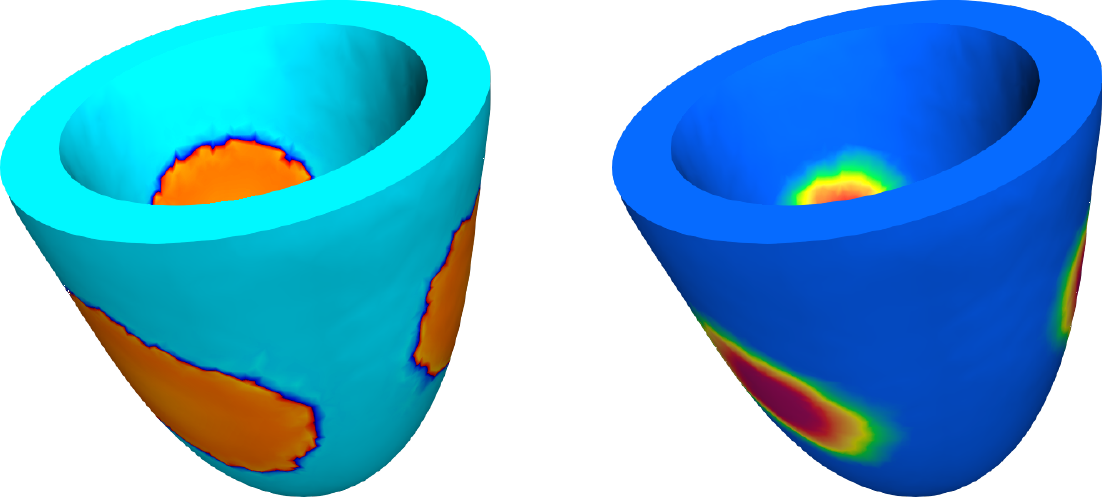}
        \caption{$t = \SI{30}{\milli\second}$}
    \end{subfigure}
    \hspace{0.07\textwidth}
    \begin{subfigure}{0.4\textwidth}
        \includegraphics[width=\textwidth]{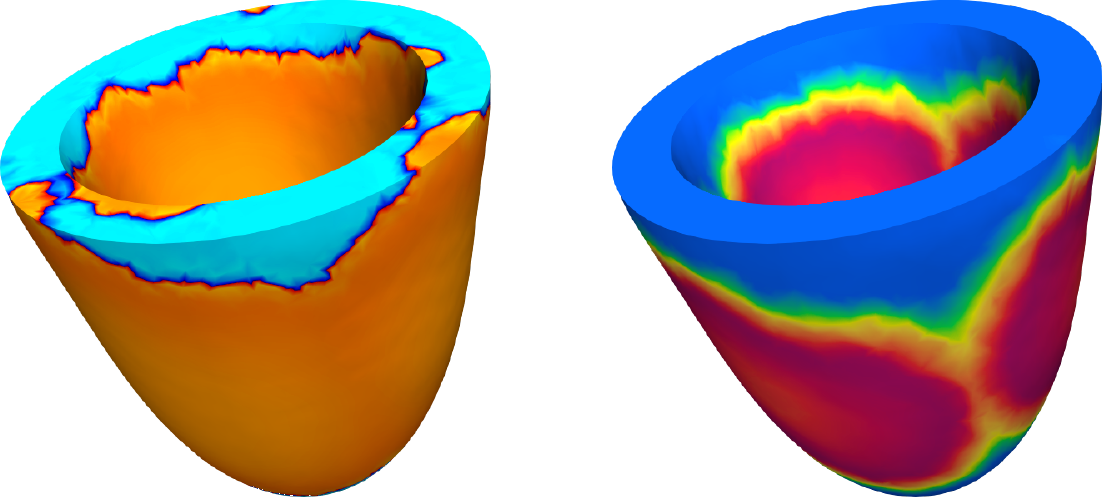}
        \caption{$t = \SI{60}{\milli\second}$}
    \end{subfigure}

    \vspace{2em}
    \begin{subfigure}{0.4\textwidth}
        \includegraphics[width=\textwidth]{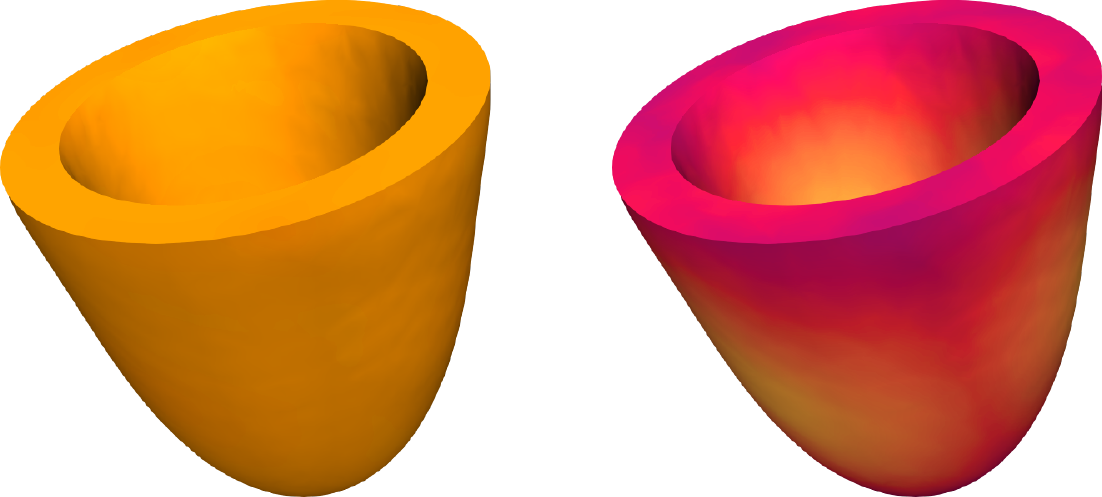}
        \caption{$t = \SI{100}{\milli\second}$}
    \end{subfigure}
    \hspace{0.07\textwidth}
    \begin{subfigure}{0.4\textwidth}
        \includegraphics[width=\textwidth]{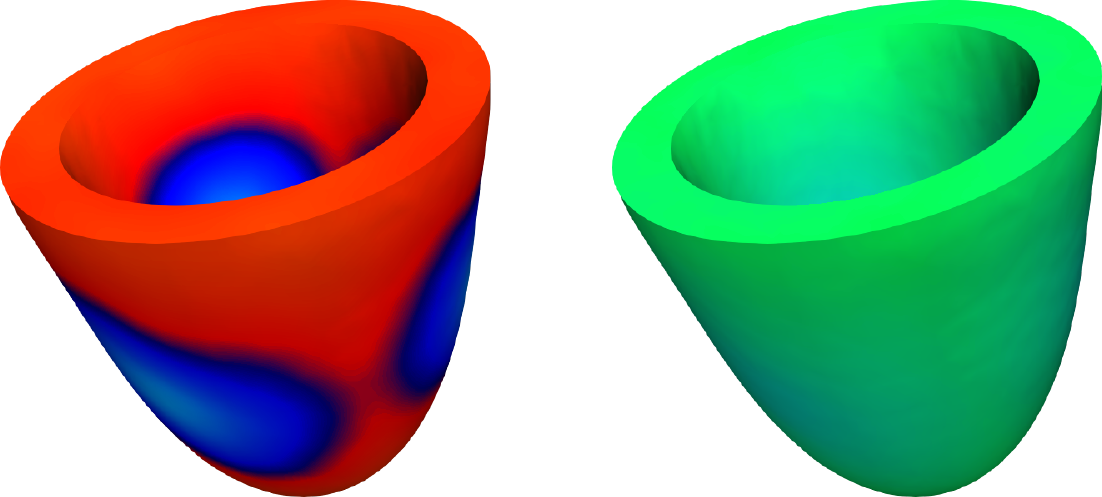}
        \caption{$t = \SI{300}{\milli\second}$}
    \end{subfigure}

    \vspace{2em}
    \begin{subfigure}{0.2\textwidth}
        \includegraphics[width=\textwidth]{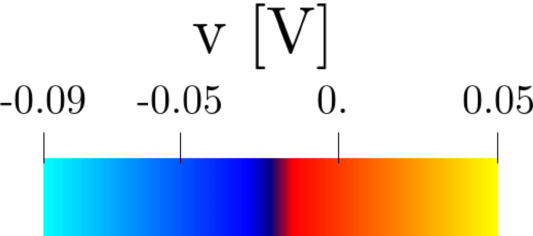}
    \end{subfigure}
    \hspace{0.05\textwidth}
    \begin{subfigure}{0.2\textwidth}
        \includegraphics[width=\textwidth]{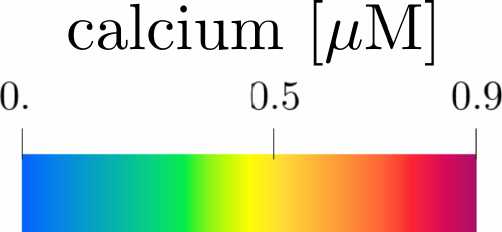}
    \end{subfigure}

    \caption{Test A. Transmembrane potential $v$ (left) and intracellular calcium concentration $[\text{Ca}^{2+}]_\text{i}$ (right) at several instants during the simulation, computed using the \exex{} with $\Delta t = \SI{0.2}{\milli\second}$.}
    \label{fig:solution-ep}
\end{figure}

\begin{figure}
    \centering

    \includegraphics[width=0.9\textwidth]{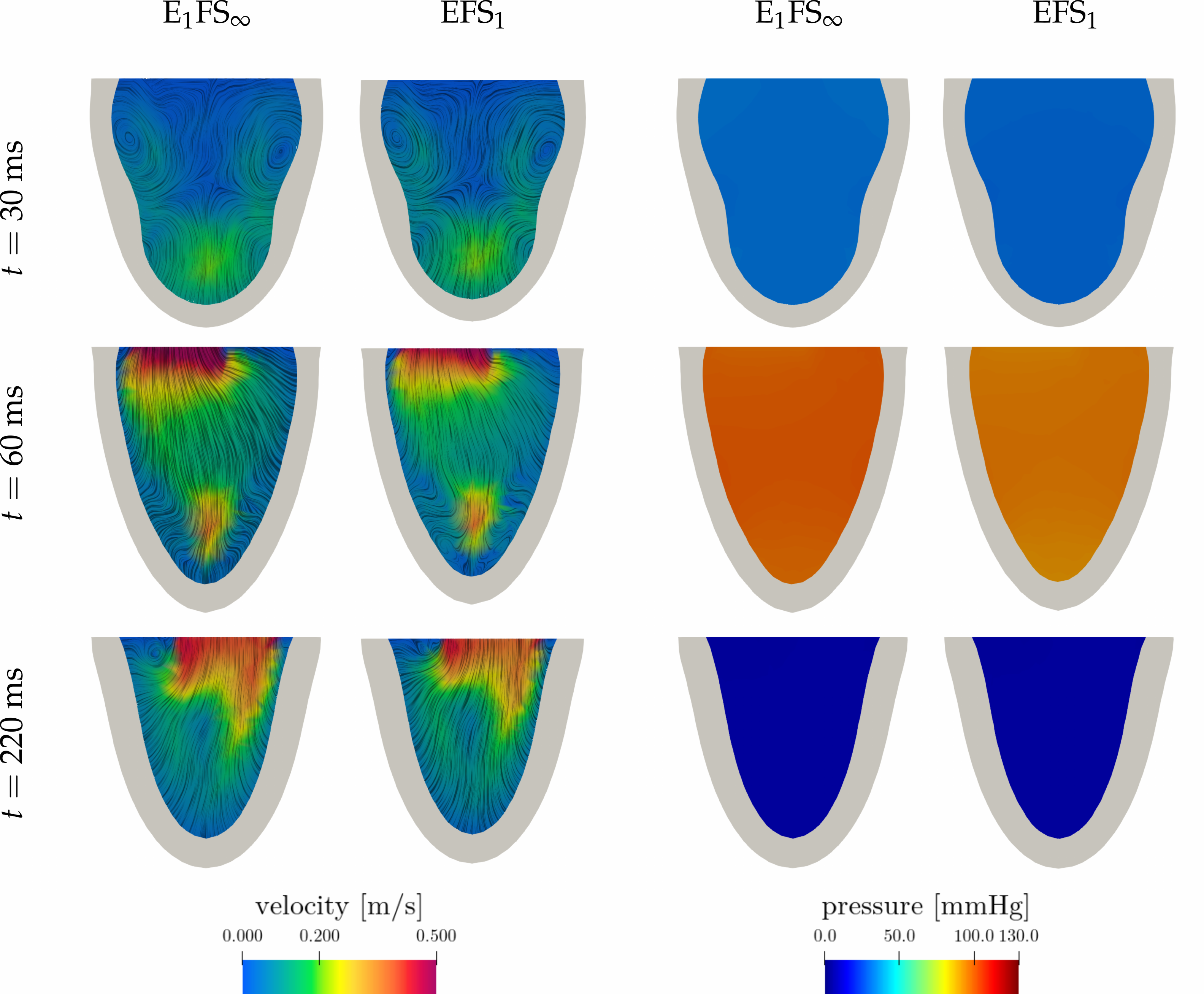}

    \caption{Test A. Fluid velocity magnitude $|\mathbf u|$ (left) and pressure $p$ (right) at three time instants computed with the \exim{} and \exex{} schemes. The velocity is overlaid with a surface line integral convolution rendering of the velocity field \cite{cabral1993imaging}. From top to bottom, the snapshots correspond are taken during the isovolumetric contraction, ejection, and filling phases, respectively.}
    \label{fig:solution-fluid}
\end{figure}

The explicit treatment of FSI coupling with RN interface conditions was shown to be conditionally stable \cite{gigante2021stability} in idealized settings. The stability condition depends on the choice of $\alpha$, on the mesh size and on the time step $\Delta t$. In order to verify numerically whether the regime and the discretization settings typical of cardiac modeling fall within the stability range, we perform tests on the idealized left ventricular geometry depicted in \cref{fig:domain} over a whole heartbeat, including all four phases. We consider the mesh $\mathcal{M}_1$ (see \cref{tab:mesh}), composed of hexahedral elements. The electrophysiology problem \eqref{eq:monodomain-discrete} is solved on a finer mesh $\mathcal{M}_1^\text{EP}$, nested into $\mathcal{M}_1$ and with half its mesh size, to better capture the sharp propagating activation front \cite{regazzoni2022cardiac,salvador2020intergrid}, and displacement and calcium are interpolated between the two meshes. We set $\Delta t = \SI{0.2}{\milli\second}$, and choose $\alpha = \SI{5000}{\sialpha}$. The value of $\alpha$ was manually tuned, starting from an initial guess derived from \cite{gigante2021choice}.

We report in \cref{fig:solution-ep,fig:solution-fluid} the solution at several time instants, computed using the \exex{} scheme. In the latter figure, the solution obtained with the \exim{} scheme is also reported. The corresponding ventricular volume and pressure over time can be found in \cref{fig:pv-alpha}. We can appreciate how in this setting the \exex{} scheme, despite being explicit, yields results that are stable in time and in qualitative agreement with those obtained with the \exim{} scheme.

In agreement with \cite{gigante2021choice,gigante2021stability}, we found the stability of the \exex{} scheme to depend on the choice of $\alpha$. Indeed, as $\alpha \to \infty$, interface conditions \eqref{eq:interface-rn} tend to Dirichlet-Neumann (DN) interface conditions \eqref{eq:interface}, which are known to lead to unstable loosely coupled schemes in the hemodynamic regime \cite{causin2005added}. As a consequence, we can expect the \exex{} method to become unstable for values of $\alpha$ not small enough. Our numerical experiments indicate that, in this setting, the \exex{} scheme is stable for all $\alpha < \SI{6750}{\sialpha}$, in qualitative accordance with \cite{gigante2021choice}.

\subsection{Test B: on the accuracy of the loosely coupled scheme}
\label{sec:test-b}

In the following sections, we consider the same setting as in Test A (\cref{sec:test-a}), and assess the accuracy of the \exex{} scheme, depending on the choice of the Robin coefficient $\alpha$ and of the time discretization step $\Delta t$.

\subsubsection{Test B1: on the influence of the Robin coefficient \texorpdfstring{$\alpha$}{} on the accuracy}

\begin{figure}
    \centering
    \begin{subfigure}{0.7\textwidth}
        \centering

        \includegraphics[width=0.49\textwidth,trim={0.5in 0in 5.8in 0in},clip]{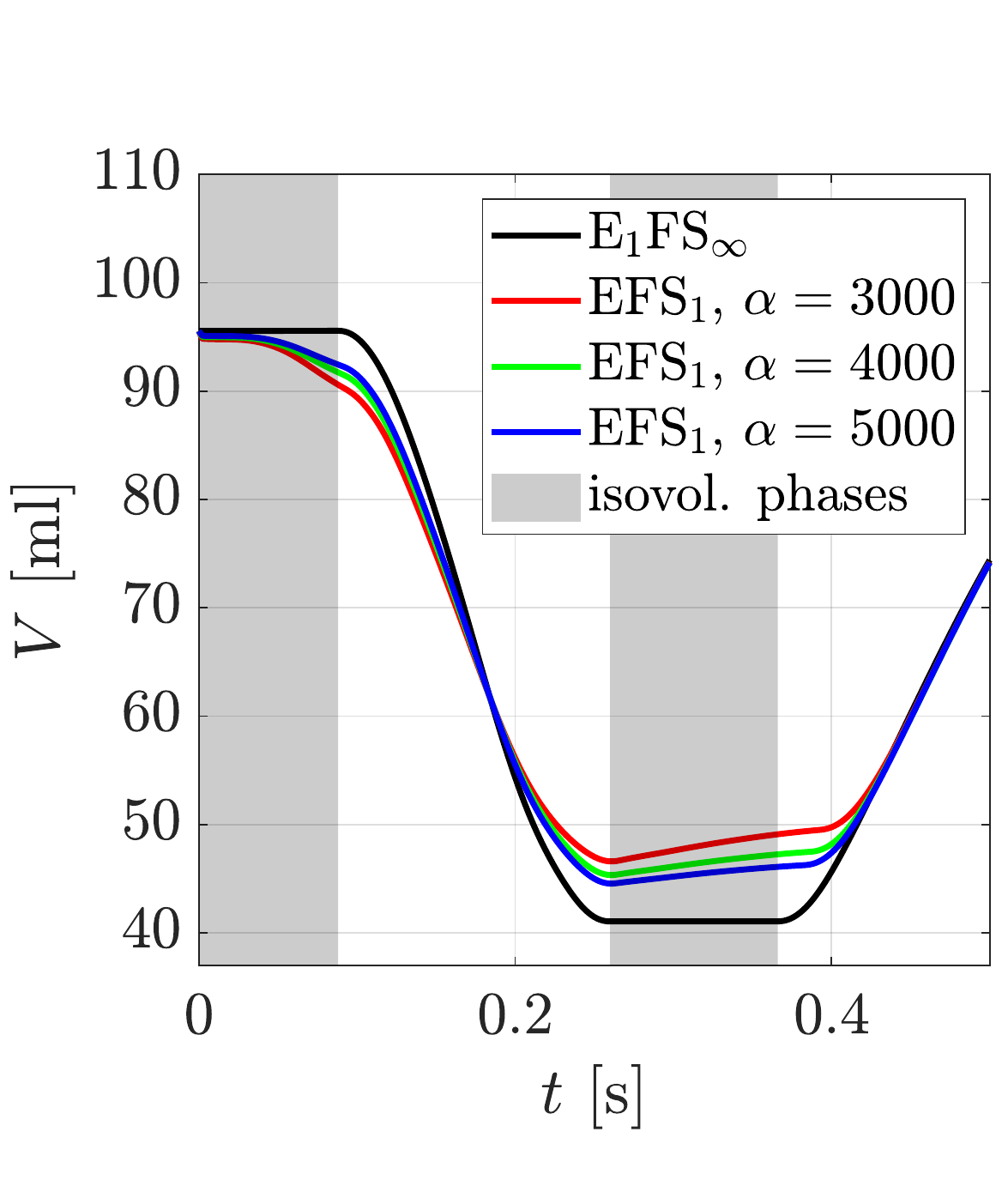}
        \includegraphics[width=0.49\textwidth,trim={5.3in 0in 1.0in 0in},clip]{fig/fig04-a.pdf}

        \caption{}
        \label{fig:pv-alpha}
    \end{subfigure}

    \begin{subfigure}{0.7\textwidth}
        \centering

        \includegraphics[width=0.49\textwidth,trim={0.5in 0in 5.8in 0in},clip]{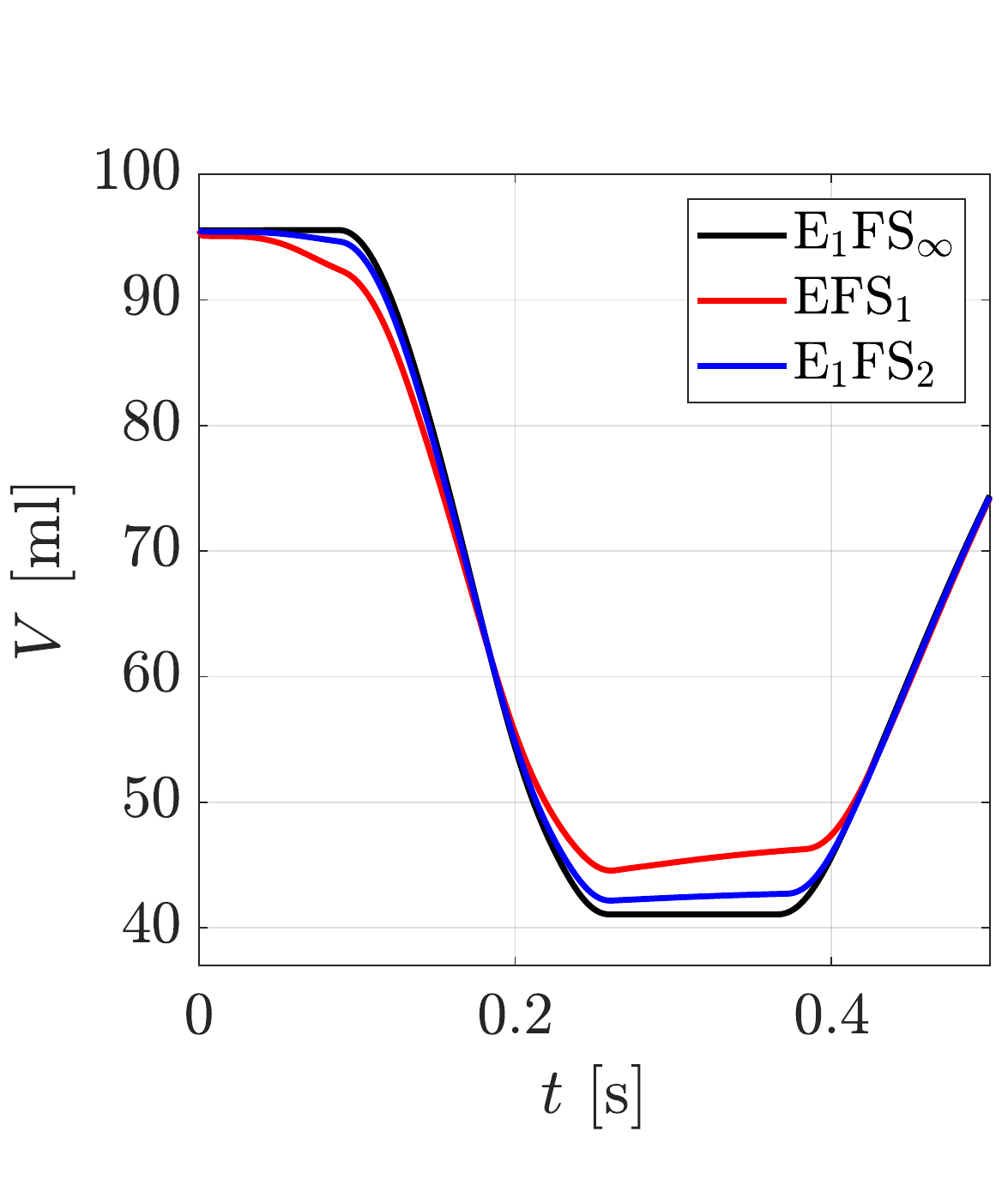}
        \includegraphics[width=0.49\textwidth,trim={5.3in 0in 1.0in 0in},clip]{fig/fig04-b.pdf}

        \caption{}
        \label{fig:rn2}
    \end{subfigure}

    \caption{(a) Test B1. Time evolution of ventricular volume (left) and average pressure (right) with the \exim{} and \exex{} schemes, with different values of the Robin coefficient $\alpha$. Grey areas identify the isovolumetric phases. (b) Test B3. Ventricular volume (left) and average pressure (right) for the \exim{}, \exin{2} and \exex{} schemes. For the last two schemes we use $\alpha = \SI{5000}{\sialpha}$.}
\end{figure}

\begin{figure}
    \centering

    \begin{subfigure}{0.49\textwidth}
        \includegraphics[width=\textwidth]{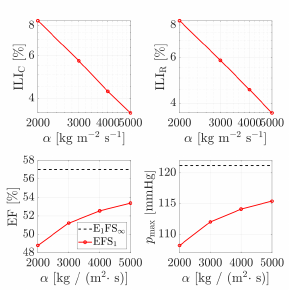}
        \caption{}
        \label{fig:indicators-vs-alpha}
    \end{subfigure}
    \begin{subfigure}{0.49\textwidth}
        \includegraphics[width=\textwidth]{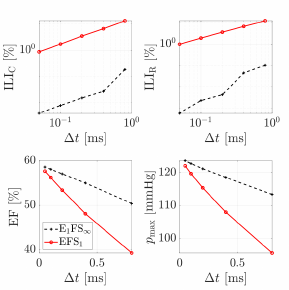}
        \caption{}
        \label{fig:indicators-vs-dt}
    \end{subfigure}

    \caption{(a) Test B1. Isovolumetric loss indices, ejection fraction, and peak systolic pressure as a function of the Robin coefficient $\alpha$. Where present, the dashed lines represent the values obtained with the explicit-implicit scheme. (b) Test B2. Isovolumetric loss indices, ejection fraction, and peak systolic pressure as a function of $\Delta t$, for the \exim{} (black) and \exex{} (red) schemes, with $\alpha = \SI{5000}{\sialpha}$.}
\end{figure}

We start by comparing the results of the \exex{} scheme against those of the \exim{} scheme, varying the Robin coefficient $\alpha$ in the range of stability experienced in Test A (i.e. $\alpha < \SI{6750}{\sialpha}$).
The choice of the parameter $\alpha$ influences the accuracy of the method. Indeed, as $\alpha \to 0$, interface conditions \eqref{eq:interface-rn} reduce to two Neumann-type conditions, and no kinematic coupling is present anymore, hindering the accuracy of the scheme.

In \cref{fig:pv-alpha} we report the time evolution of ventricular volumes and pressures for the \exim{} scheme and for different values of $\alpha$ in the \exex{} scheme. From these results, we can observe a general qualitative agreement between \exex{} and \exim{} solutions. However, unlike the latter, the \exex{} solutions feature a loss of mass and, consequently, volume variations during the isovolumetric phase. This leads to a slower ejection and filling, as well as a lower peak pressure. This behavior increases for decreasing values of $\alpha$.

Similar conclusions can be drawn by looking at the plots in \cref{fig:indicators-vs-alpha}, where the value of the considered indicators obtained by the \exex{} scheme has been plotted against the Robin coefficient $\alpha$. Although the mismatch reduces as $\alpha$ increases, even with the highest value of $\alpha$ the two results present differences of \SI{6.3}{\percent} in ejection fraction and \SI{4.7}{\percent} in peak systolic pressure.

\subsubsection{Test B2: on the influence of \texorpdfstring{$\Delta t$}{the time step} on the accuracy}

We expect the mismatch between the \exim{} and \exex{} schemes to reduce as $\Delta t$ is reduced. To this end, considering $\alpha = \SI{5000}{\sialpha}$, we perform several simulations reducing the time step of both the schemes. The resulting indicators are reported in \cref{fig:indicators-vs-dt}. As expected, we observe that as $\Delta t \to 0$ there is increasing agreement between the solutions computed by the two schemes in terms of EF and $p_{\max}$. For both schemes the isovolumetric loss indices $\text{ILI}_\text{C}$ and $\text{ILI}_\text{R}$ tend to zero as $\Delta t \to 0$, with similar rates. However, the ones obtained with the \exim{} scheme are in any case smaller than those of the \exex{} one.

\begin{figure}
    \centering

    \includegraphics[width=\textwidth]{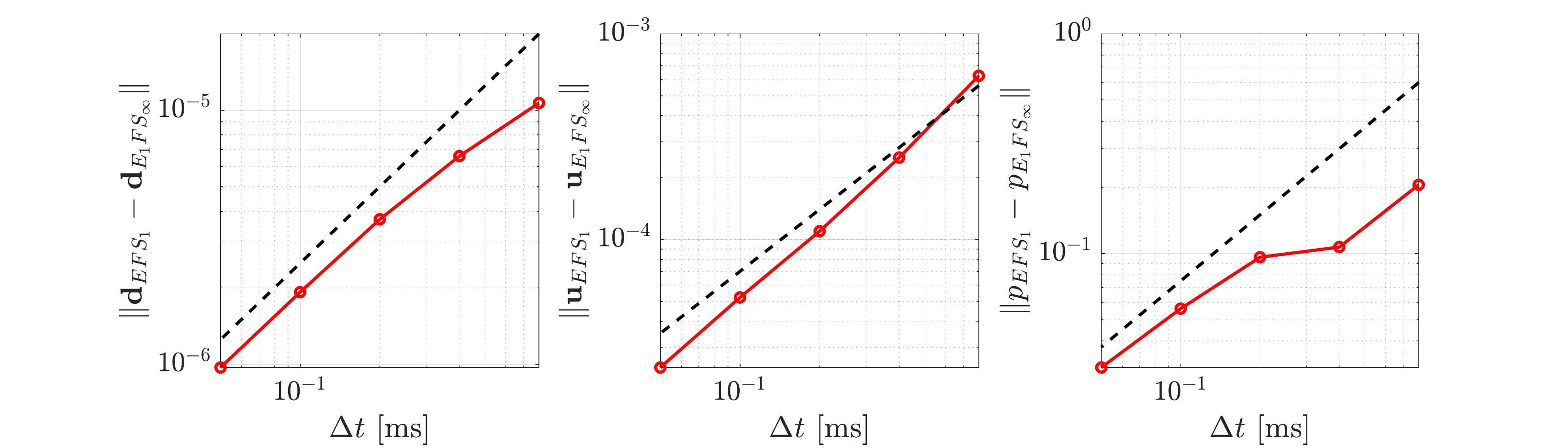}

    \caption{Norm of the difference between the solutions (from left to right, displacement, velocity and pressure) of the \exex{} and \exim{} schemes, with varying $\Delta t$. The dashed black lines are parallel to $f(\Delta t) = \Delta t$, and are used as a reference for the convergence order.}
    \label{fig:dt-convergence}
\end{figure}

We report in \cref{fig:dt-convergence} the norm of the difference between the solutions computed with the \exex{} and \exim{} schemes, with varying $\Delta t$. We observe that the mismatch tends to zero as $\Delta t \to 0$, with order \num{1}. Therefore, we conclude that the segregation of the fluid and solid solver introduces a splitting error which is at most of order \num{1}, the same order of the time discretization used for the individual subproblems. We remark that higher-order time discretization schemes may require to enhance the \exex{} scheme to preserve the time convergence order.

\subsubsection{Test B3: on the influence of multiple Robin-Neumann subiterations on the accuracy}

\begin{table}
    \centering
    \begin{tabular}{c c c c c}
        \toprule
        \textbf{Scheme} & $\text{ILI}_\text{C}$ [\%] & $\text{ILI}_\text{R}$ [\%] & $\text{EF}$ [\%] & $p_{\max}$ [\si{\mmhg}] \\
        \midrule
        \exim{} & \num{0.02} & \num{0.00} & \num{57.0} & \num{121.2} \\
        \exex{} & \num{3.48} & \num{3.66} & \num{53.4} & \num{115.4} \\
        \exin{2} & \num{0.96} & \num{1.29} & \num{55.9} & \num{119.5} \\
        \bottomrule
    \end{tabular}
    \caption{Test B3. Isovolumetric loss indices, ejection fraction and peak systolic pressure for three simulations, using the \exim{}, \exex{} and \exin{2} schemes. For the last two schemes, we set $\alpha = \SI{5000}{\sialpha}$.}
    \label{tab:rn2}
\end{table}

We run a simulation using the scheme \exin{2}, with $\alpha = \SI{5000}{\sialpha}$ and $\Delta t = \SI{0.2}{\milli\second}$, and compare the results against the corresponding ones obtained with the \exim{} and \exex{} schemes in terms of the evolutions of ventricular volume and pressure. The results obtained are reported in \cref{fig:rn2}, while in \cref{tab:rn2} we report the values of the ILI, EF and $p_{\max}$ indicators for the three schemes. We observe that doing two RN iterations can significantly improve the agreement with the explicit-implicit scheme. In particular, the isovolumetric phases are captured more accurately, as indicated by the ILI indices.

\subsubsection{Test B4: on the accuracy during the ejection phase}

\begin{table}
    \centering
    \begin{tabular}{c c c}
        \toprule
        \textbf{Scheme} & $\text{EF}$ [\%] & $p_{\max}$ [\si{\mmhg}] \\
        \midrule
        \exim{} & \num{57.0} & \num{121.2} \\
        \exex{} & \num{54.4} & \num{119.6} \\
        \exex{} (full heartbeat) & \num{53.4} & \num{115.4} \\
        \bottomrule
    \end{tabular}
    \caption{Test B4. Ejection fraction and peak systolic pressure in the ejection phase test, for the \exim{} and \exex{} schemes. For comparison, we report the same quantities computed in the full heartbeat test (Test A) with the \exex{} scheme. In both cases, we set $\alpha = \SI{5000}{\sialpha}$.}
    \label{tab:indicators-systole}
\end{table}

\begin{figure}
    \centering
    \begin{subfigure}{0.35\textwidth}
        \includegraphics[width=\textwidth,trim={0.5in 0in 5.8in 0in},clip]{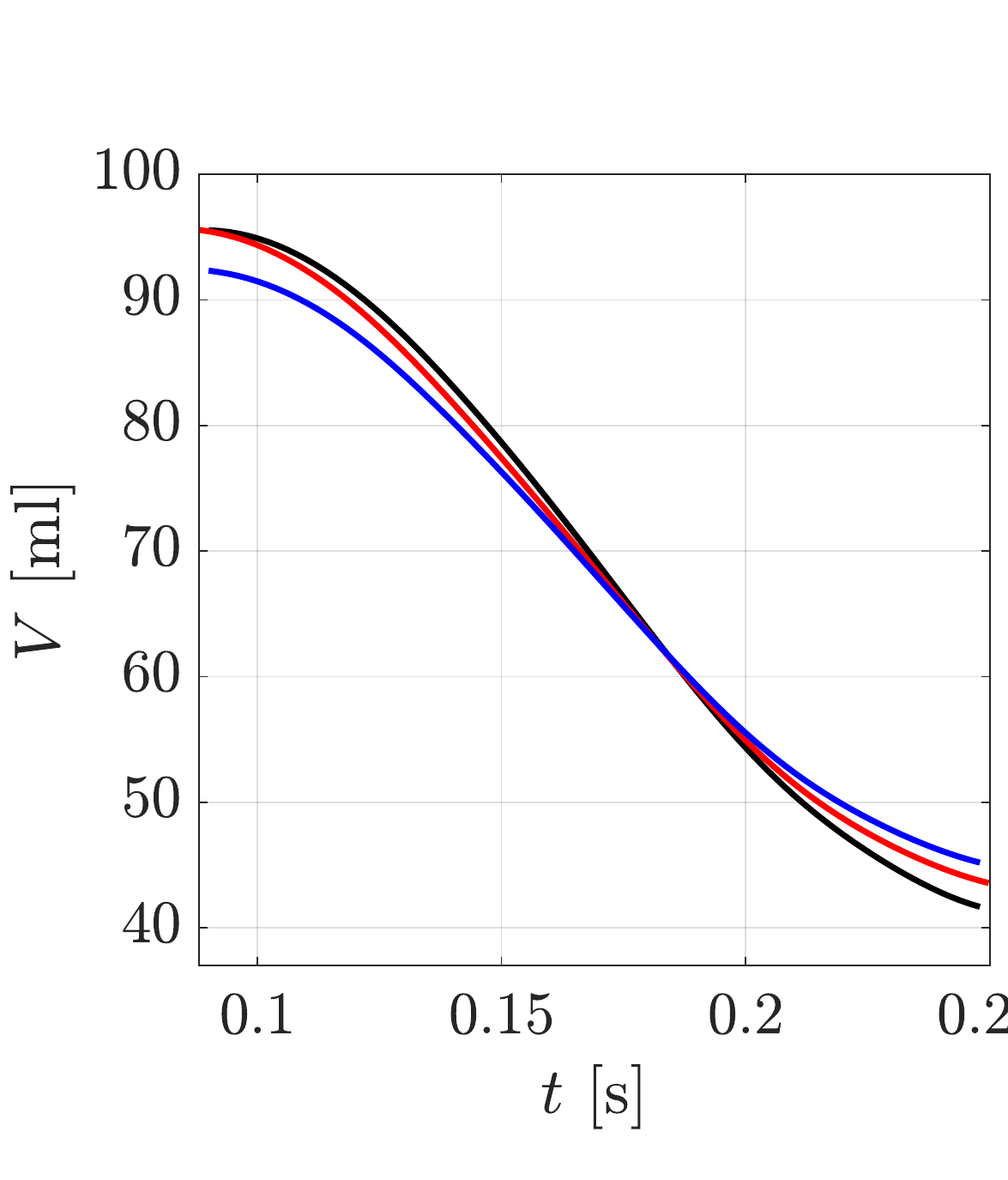}
        \caption{}
    \end{subfigure}
    \hspace{0.03\textwidth}
    \begin{subfigure}{0.35\textwidth}
        \includegraphics[width=\textwidth,trim={5.4in 0in 0.9in 0in},clip]{fig/fig07.pdf}
        \caption{}
    \end{subfigure}

    \caption{Test B4. Ventricular volume (a) and average pressure (b) in the ejection phase test, computed with the \exim{} and \exex{} schemes. For comparison, we report in blue the volume and pressure corresponding to the full heartbeat test (Test B1) with the \exex{} scheme. For loosely coupled schemes we set $\alpha = \SI{5000}{\sialpha}$.}
    \label{fig:pv-systole}
\end{figure}

Previous sections show that the \exex{} scheme introduces an error in capturing volume conservation during isovolumetric phases. This has an impact on the evolution of pressure during those phases, that in turn influences the solution during ejection and filling.

To understand to what extent the mismatch between \exim{} and \exex{} schemes is determined by isovolumetric phases, we simulate only the ejection phase, by providing as initial condition the solution of the \exim{} scheme from test A (\cref{sec:test-a}) at time $t_0 = \SI{88}{\milli\second}$ (corresponding to the end of isovolumetric contraction). We compare EF and $p_{\max}$ obtained with the \exim{} and \exex{} schemes. Results for these indicators are reported in \cref{tab:indicators-systole}. We observe that, while both ejection fraction and peak pressure are smaller in the \exex{} case than they are in the \exim{} case, the reduction is less significant than what is observed in a full heartbeat explicit simulation. Similar conclusions are drawn by comparing the pressure and volume over time, as reported in \cref{fig:pv-systole}: in the ejection-only simulation, there is better agreement between the \exim{} and the \exex{} schemes.

Overall, this result indicates that the mismatch between the two schemes could be particularly relevant during the isovolumetric phases. Therefore, it can be of interest to explore adaptive methods that adjust e.g. the number of RN subiterations depending on the simulated heartbeat phase.

\subsection{Test C: computational efficiency}

\begin{figure}
    \centering
    \begin{subfigure}{0.29\textwidth}
        \includegraphics[width=\textwidth,trim={1.0in 0in 9in 0in},clip]{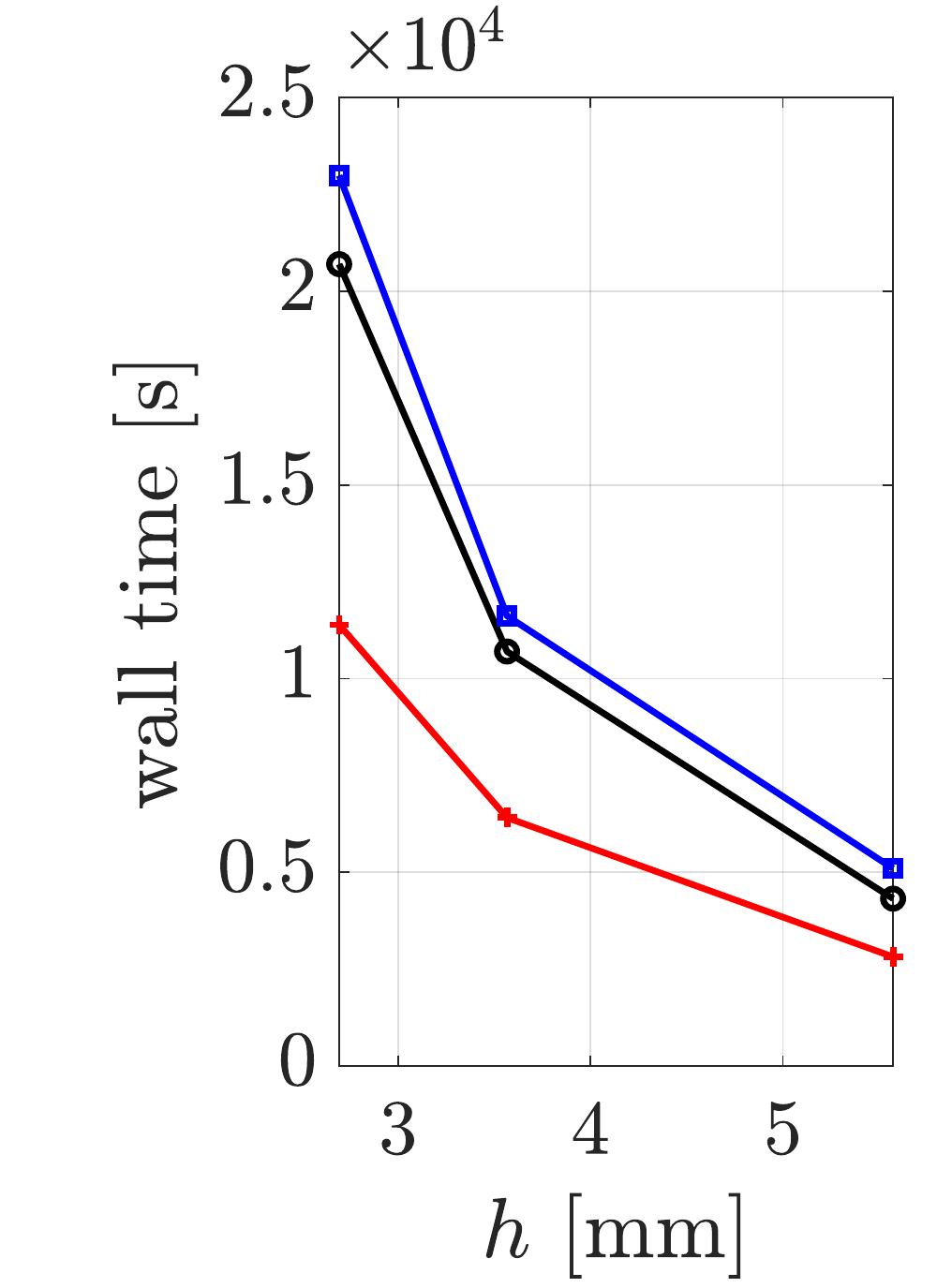}
        \caption{}
    \end{subfigure}
    \hspace{0.01\textwidth}
    \begin{subfigure}{0.29\textwidth}
        \includegraphics[width=\textwidth,trim={5.0in 0in 5.0in 0in},clip]{fig/fig08.pdf}
        \caption{}
    \end{subfigure}
    \hspace{0.01\textwidth}
    \begin{subfigure}{0.29\textwidth}
        \includegraphics[width=\textwidth,trim={9.0in 0in 1.0in 0in},clip]{fig/fig08.pdf}
        \caption{}
    \end{subfigure}

    \caption{Test C. Total wall time (a), wall time spent to assemble the fluid and structure systems (b) and to solve them (c) against average mesh size $h$ of the three considered meshes (\cref{tab:mesh}). \exim{}, \exex{} and \exin{2} schemes. For the last two schemes, we set $\alpha = \SI{5000}{\sialpha}$.}
    \label{fig:walltime-vs-h}
\end{figure}

The chief advantage of a loosely coupled scheme is its computational efficiency if compared to a method where the couplings are treated implicitly. To verify this, we perform numerical simulations of the \exim{}, \exex{} and \exin{2} schemes with three differently refined meshes, detailed in \cref{tab:mesh}. We compare the total wall time, the portion of wall time devoted to the assembly of fluid and structure systems, as well as the wall time spent in the solution of the fluid, structure or FSI systems. We do not consider in detail computational times associated to electrophysiology, force generation and fluid domain displacement, since the considered schemes are identical in those steps. These simulations ran in parallel using \num{44} cores with Intel Xeon Platinum 8160@2.1GHz processors.

Results are reported in \cref{fig:walltime-vs-h}. From these results, we can appreciate how the \exex{} scheme leads to a very significant reduction in computational time with respect to the \exim{} one. This reduction becomes more significant as the mesh is refined: the total wall time for the simulation on the finest mesh $\mathcal{M}_3$ using the \exex{} scheme is approximately \SI{45}{\percent} less than the corresponding simulation using the \exim{} scheme. In particular, the cost associated to both the assembly and the solution of the linear systems for the FSI problem is much smaller in the \exex{} scheme than it is in the \exim{}. Overall, the \exex{} scheme allows for a significant saving in computational time with respect to the \exim{} one. Conversely, the \exin{2} requires a computational time similar to that of the \exim{} scheme.

\subsection{Test D: the case of a realistic human ventricle}
\label{sec:zygote}

\begin{figure}
    \centering
    \begin{subfigure}{0.30\textwidth}
        \includegraphics[width=\textwidth]{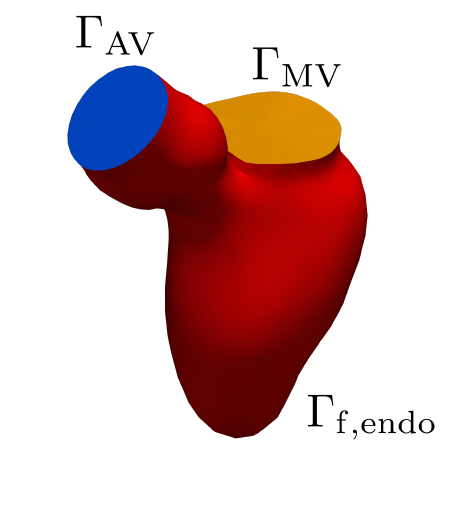}
        \caption{}
    \end{subfigure}
    \begin{subfigure}{0.30\textwidth}
        \includegraphics[width=\textwidth]{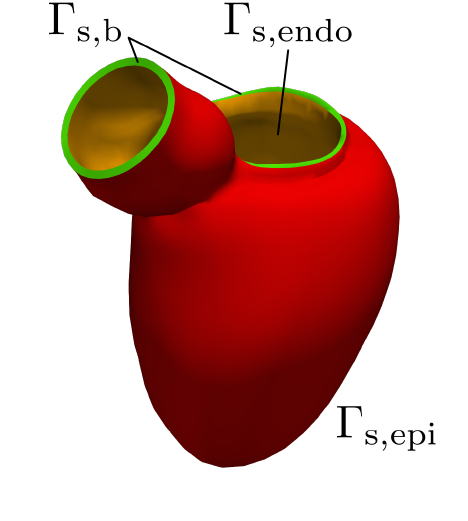}
        \caption{}
    \end{subfigure}
    \begin{subfigure}{0.30\textwidth}
        \includegraphics[width=\textwidth]{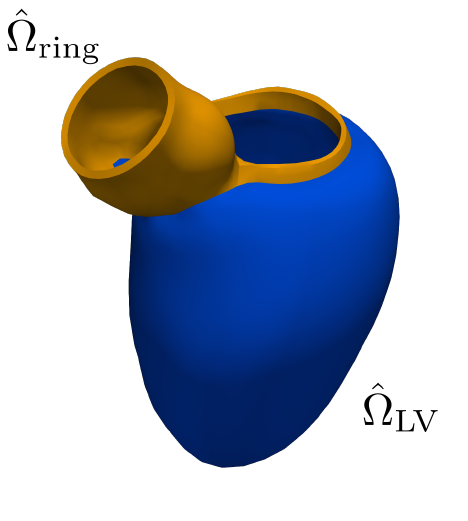}
        \caption{}
        \label{fig:domain-zygote-subdomains}
    \end{subfigure}
    \caption{Test D. Computational domain $\Omega\fluid$ (a) and $\Omega\solid$ (b, c) of the realistic ventricle. Colors and labels denote the different portions of the boundary in (a) and (b), and denote the two subdomains in (c).}
    \label{fig:zygote-domain}
\end{figure}

We present a test case in a more realistic setting to showcase the effectiveness of the proposed scheme. We consider the left ventricle from the heart model provided by Zygote Media Group \cite{zygote}, represented in \cref{fig:zygote-domain}. We processed the geometry using the meshing algorithms presented in \cite{fedele2021polygonal} using the software \texttt{VMTK} \cite{vmtk}.

The model includes ventricular inflow and outflow tracts. Those portions are not formed of muscular tissue, as the bulk of the ventricle is \cite{katz2010physiology}. To account for this, we introduce two subdomains into $\hat\Omega\solid$, denoted by $\hat\Omega_\text{LV}$ and $\hat\Omega_\text{ring}$ (see \cref{fig:domain-zygote-subdomains}), representing the left ventricle and the valvular rings, respectively, and employ a neo-Hookean constitutive law in $\hat\Omega_\text{ring}$ (while keeping the Guccione constitutive law in $\hat\Omega_\text{LV}$). We also set $T_\text{act,max} = 0$ in $\hat\Omega_\text{ring}$.

For the discretization, we use tetrahedral elements for fluid and solid domains. The mechanics and fluid dynamics equations are discretized with linear finite elements. To deal with the higher accuracy requirements of electrophysiology, we use quadratic finite elements to discretize \cref{eq:monodomain-discrete}. This is an alternative approach to the one used in previous sections, based on nested mesh refinement. We set $\Delta t = \SI{0.2}{\milli\second}$. The test ran on 48 cores from the CINECA GALILEO100 supercomputer\footnote{Technical specifications: \url{https://wiki.u-gov.it/confluence/display/SCAIUS/UG3.3\%3A+GALILEO100+UserGuide}}.

\begin{figure}
    \centering
    \begin{subfigure}{0.4\textwidth}
        \includegraphics[width=\textwidth]{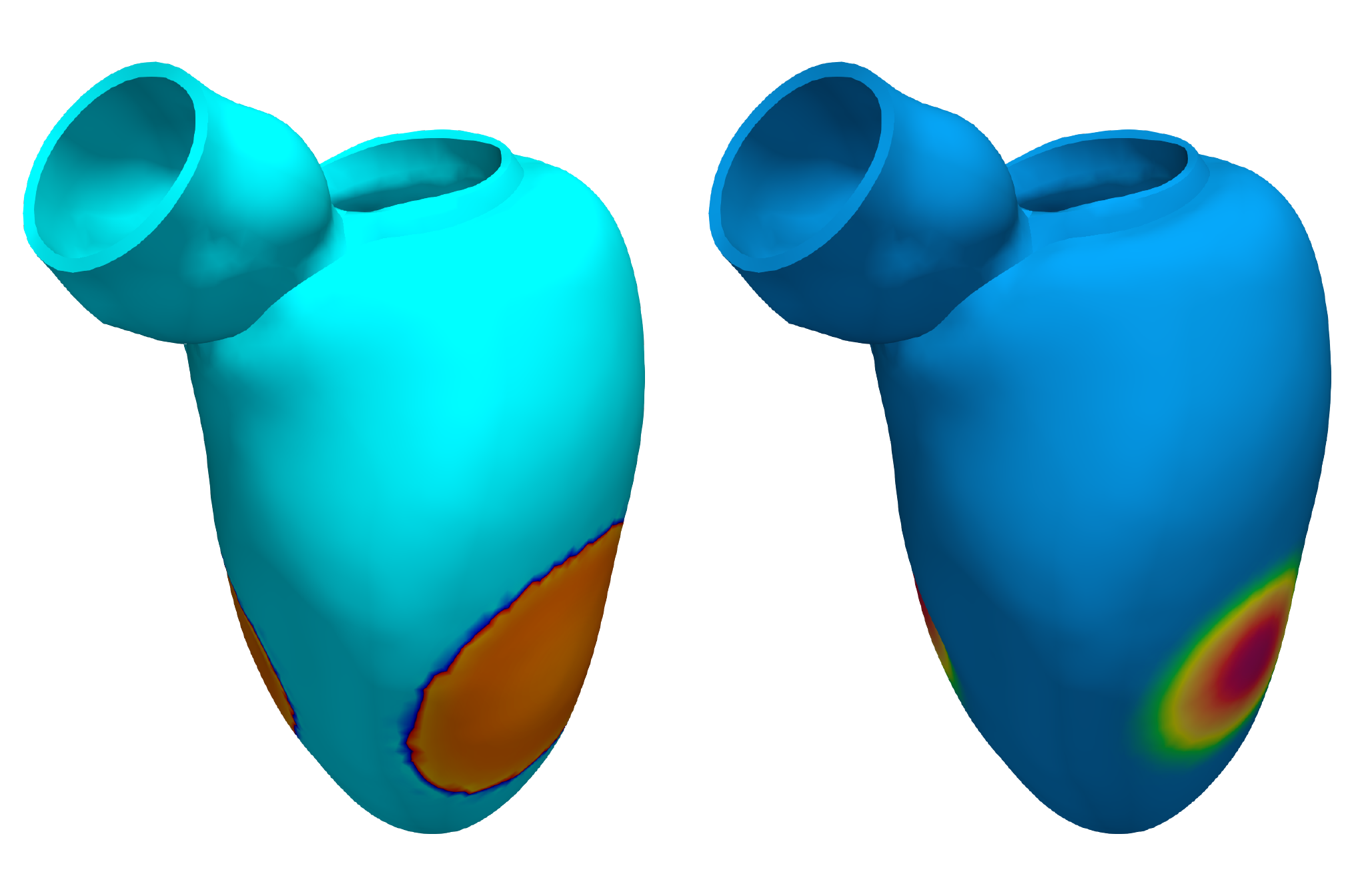}
        \caption{$t = \SI{30}{\milli\second}$}
    \end{subfigure}
    \hspace{0.022\textwidth}
    \begin{subfigure}{0.4\textwidth}
        \includegraphics[width=\textwidth]{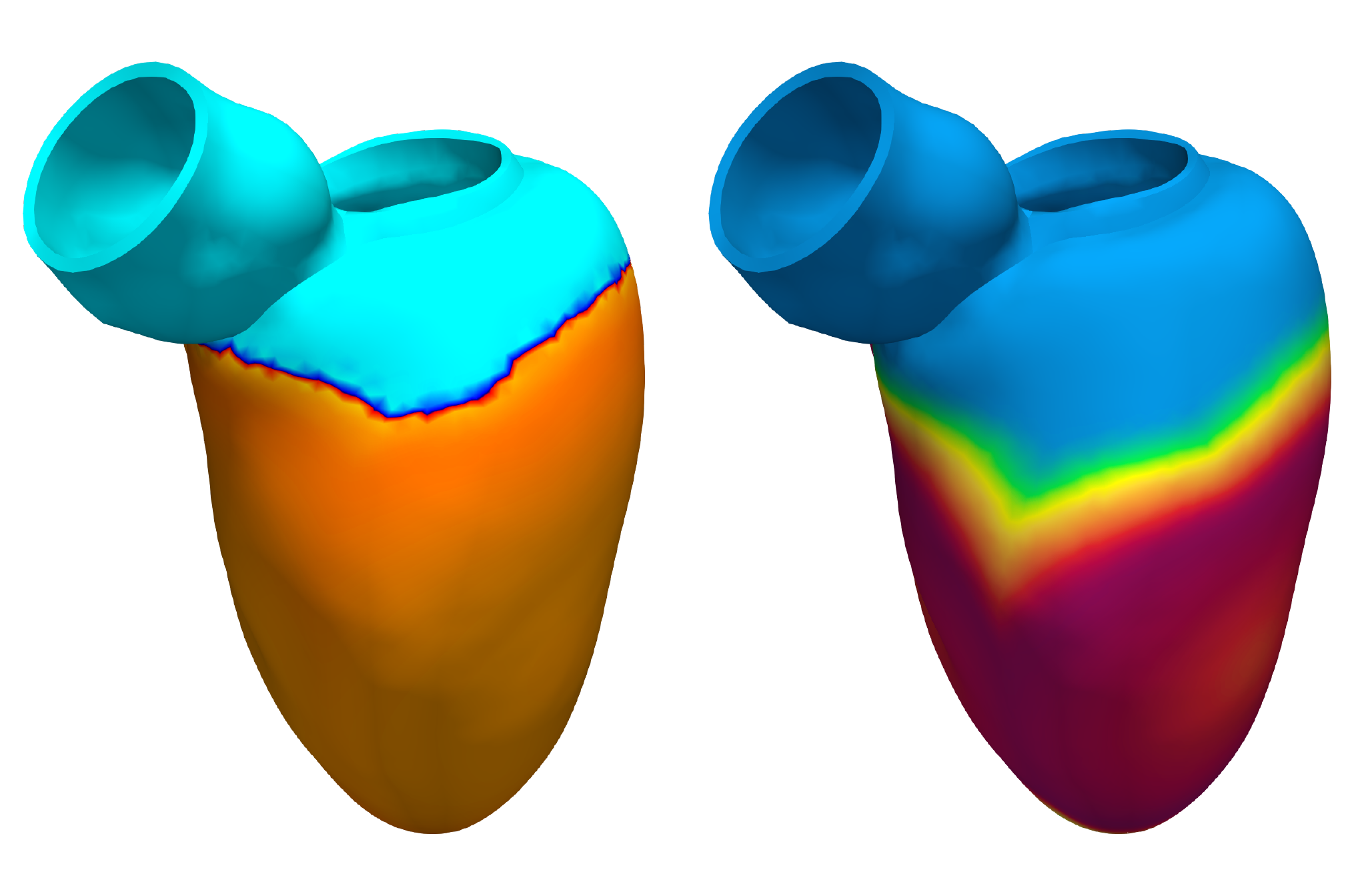}
        \caption{$t = \SI{70}{\milli\second}$}
    \end{subfigure}
    \hspace{0.022\textwidth}
    \begin{subfigure}{0.4\textwidth}
        \includegraphics[width=\textwidth]{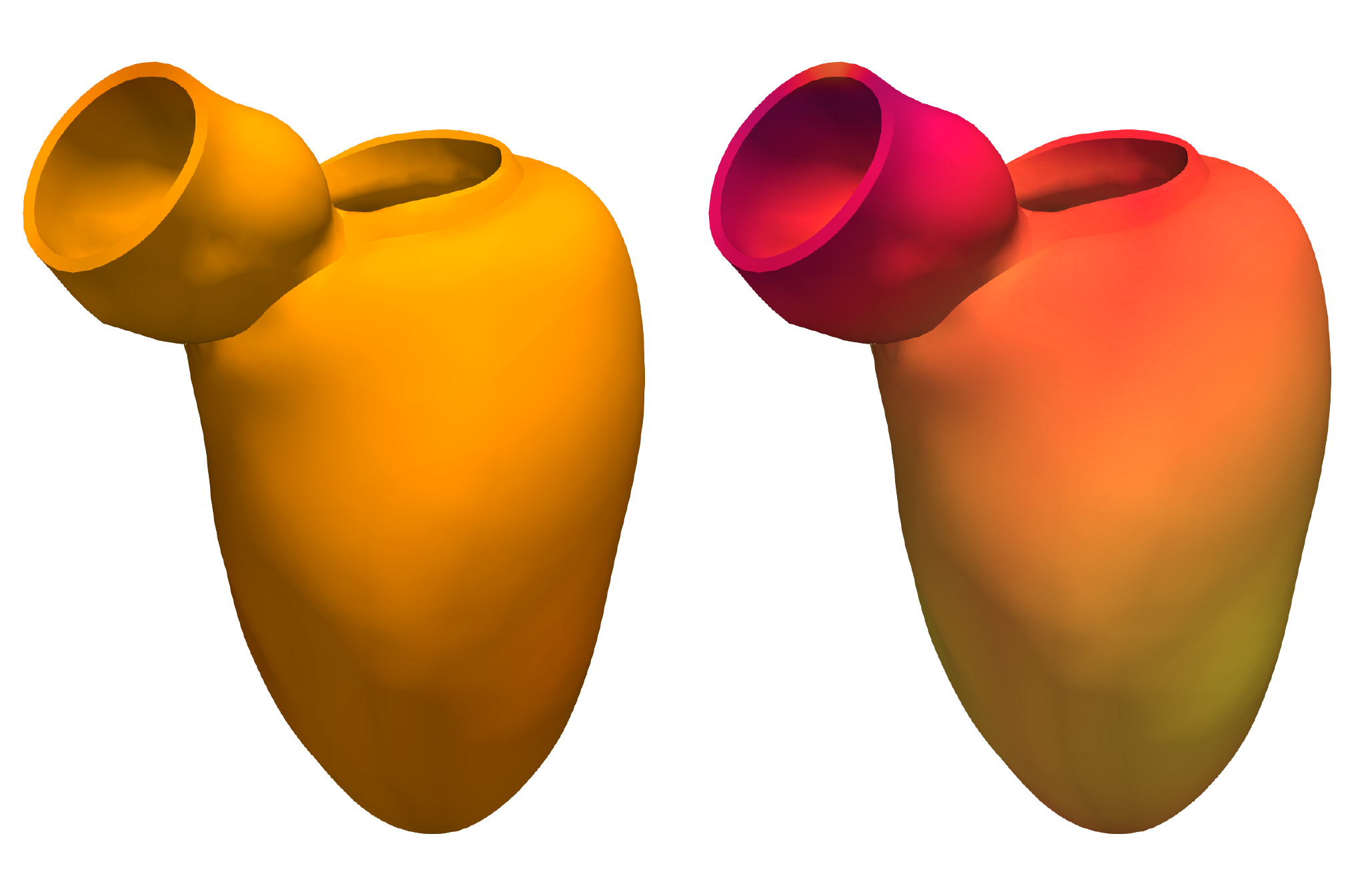}
        \caption{$t = \SI{150}{\milli\second}$}
    \end{subfigure}
    \hspace{0.022\textwidth}
    \begin{subfigure}{0.4\textwidth}
        \includegraphics[width=\textwidth]{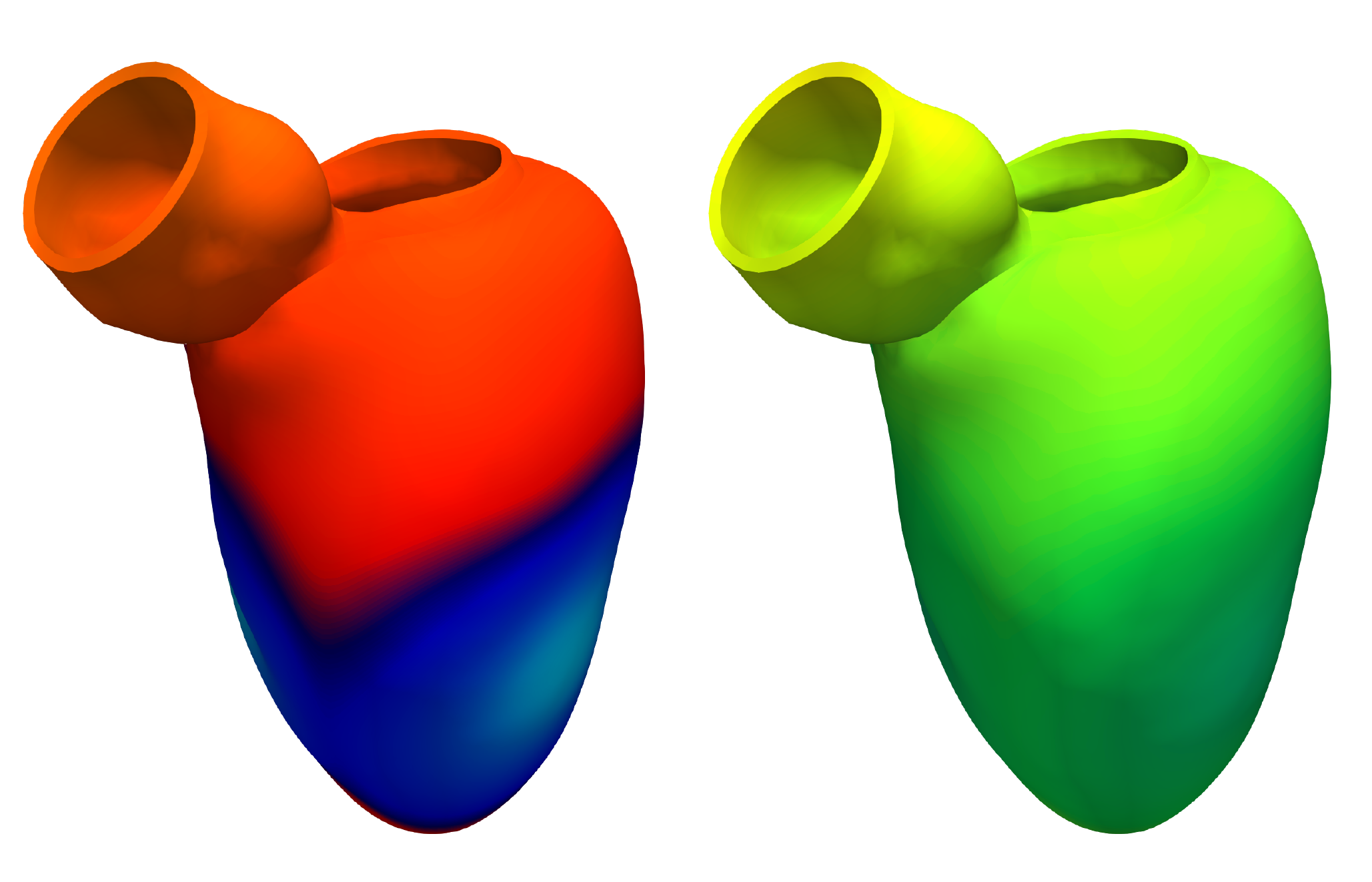}
        \caption{$t = \SI{310}{\milli\second}$}
    \end{subfigure}

    \begin{subfigure}{0.2\textwidth}
        \includegraphics[width=\textwidth]{fig/fig02-e.pdf}
    \end{subfigure}
    \hspace{0.05\textwidth}
    \begin{subfigure}{0.2\textwidth}
        \includegraphics[width=\textwidth]{fig/fig02-f.pdf}
    \end{subfigure}

    \caption{Test D. Transmembrane potential $v$ (left) and intracellular calcium concentration $[\text{Ca}^{2+}]_\text{i}$ (right) at several instants during the simulation of the realistic ventricle, computed using the \exex{} scheme.}
    \label{fig:solution-zygote-ep}
\end{figure}

\begin{figure}
    \centering
    \includegraphics[width=0.8\textwidth]{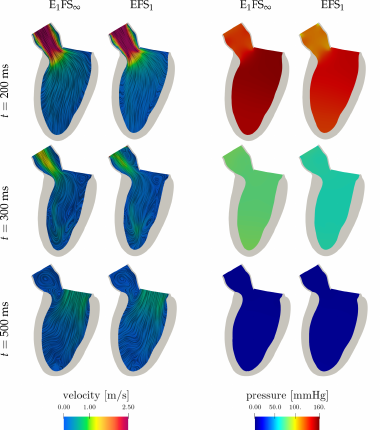}

    \caption{Test D. Fluid velocity magnitude $|\mathbf u|$ (left) and pressure $p$ (right) at three instants during the simulation of the realistic ventricle, computed using the \exim{} and \exex{} schemes. The velocity magnitude is overlaid with a surface line integral convolution rendering of the flow field \cite{cabral1993imaging}.}
    \label{fig:solution-zygote-fluid}
\end{figure}

\begin{figure}
    \centering
    \begin{subfigure}{0.35\textwidth}
        \includegraphics[width=\textwidth,trim={0.5in 0in 5.8in 0in},clip]{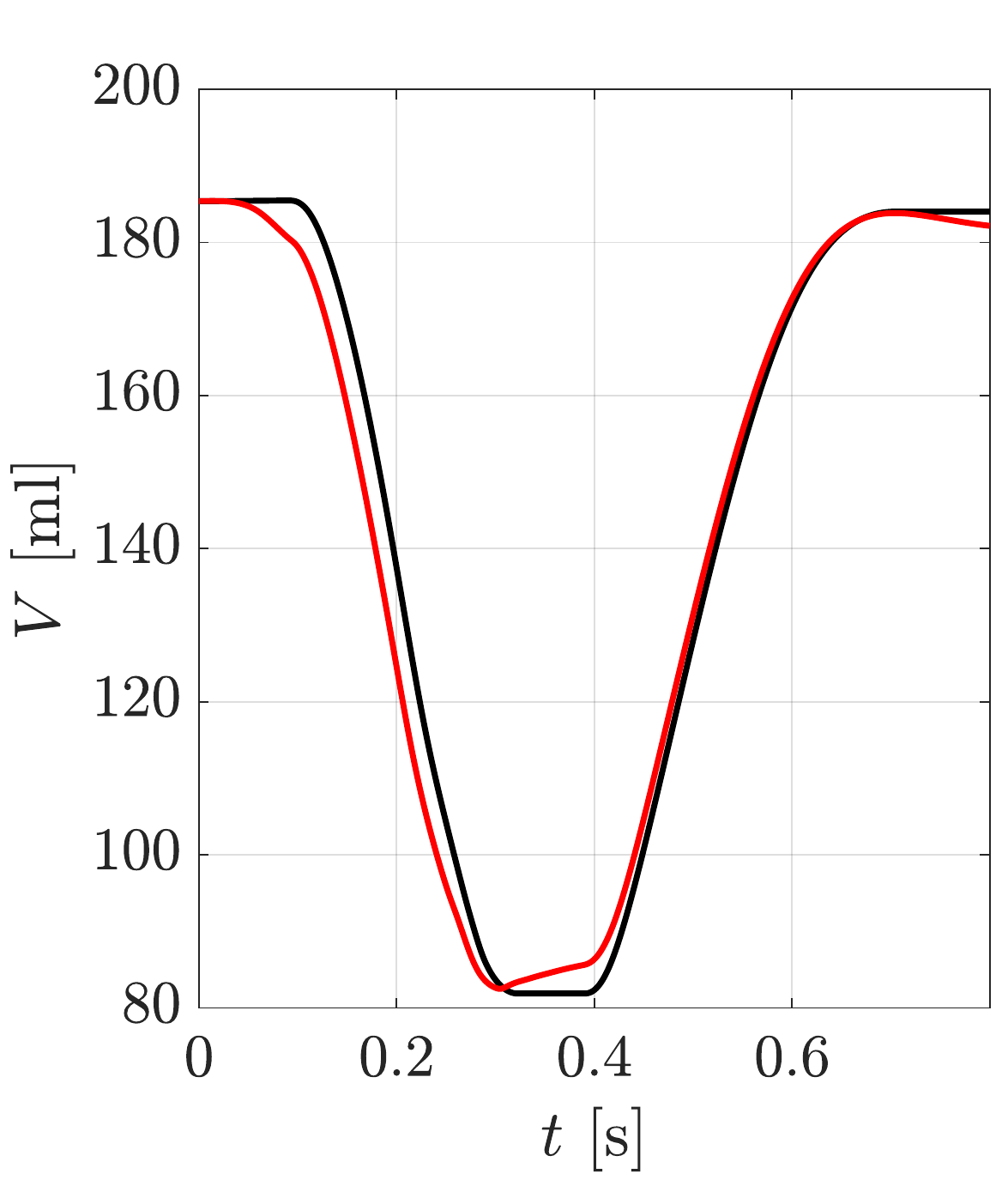}
        \caption{}
    \end{subfigure}
    \hspace{0.03\textwidth}
    \begin{subfigure}{0.35\textwidth}
        \includegraphics[width=\textwidth,trim={5.3in 0in 1.0in 0in},clip]{fig/fig12.pdf}
        \caption{}
    \end{subfigure}

    \caption{Test D. Ventricular volume (a) and average pressure (b) in the realistic test for the \exim{} and \exex{} schemes, with $\alpha = \SI{5000}{\sialpha}$.}
    \label{fig:pv-zygote}
\end{figure}

We report in \cref{fig:solution-zygote-ep} some snapshots of the electrophysiology solution for this test case, while in \cref{fig:solution-zygote-fluid} we report a comparison of domain deformation and fluid dynamics variables with the \exim{} solution. From these results we observe the stability of the proposed loosely coupled scheme and the qualitative agreement of the solution with the \exim{} one.

In \cref{fig:pv-zygote} we show the ventricular volume and pressure over time for the two schemes, whereas in \cref{tab:zygote-indicators} we report the computed indicators. We observe again that the \exex{} scheme introduces an error in capturing the isovolumetric phases, and this leads to a mismatch with the \exim{} scheme in terms of peak systolic pressure, while the ejection fraction is well captured.

\begin{table}
    \centering
     \begin{tabular}{c c c c c}
        \toprule
        \textbf{Scheme} & $\text{ILI}_\text{C}$ [\%] & $\text{ILI}_\text{R}$ [\%] & $\text{EF}$ [\%] & $p_{\max}$ [\si{\mmhg}] \\
        \midrule
        \exim{} & \num{0.04} & \num{0.00} & \num{55.8} & \num{148.9} \\
        \exex{} & \num{2.76} & \num{3.60} & \num{55.5} & \num{136.2} \\
        \bottomrule
    \end{tabular}
    \caption{Test D. Isovolumetric loss indicators, ejection fraction and peak systolic pressure for the realistic test case, using the \exim{} and \exex{} schemes, with $\alpha = \SI{5000}{\sialpha}$.}
    \label{tab:zygote-indicators}
\end{table}

\section{Conclusions}
\label{sec:conclusions}

We propose a loosely coupled scheme (\exex{}) in the context of cardiac simulations, for the coupling of electrophysiology, active and passive tissue mechanics, and hemodynamics, where the subproblems are solved only once per time step and a Robin interface condition is considered for the fluid subproblem. We compare its performance with two other schemes where electrophysiology is treated explicitly, one with strong FSI coupling (\exim{}) and one where \num{2} FSI iterations are performed (\exin{2}). The main findings of our work (valid for both an idealized and a realistic geometry) are:
\begin{enumerate}
    \item the \exex{} scheme is stable in the physiological regime, provided that the Robin interface parameter $\alpha$ is small enough;
    \item \exex{} introduces a further error, besides the time discretization, due to the lack of synchrony between the FSI interface conditions. This error is mostly relevant during the isovolumetric phases, although it vanishes for decreasing values of the time step $\Delta t$;
    \item \exex{} is about \SI{45}{\percent} faster than \exim{}, with the computational saving becoming more relevant as the mesh is refined.
\end{enumerate}
In conclusion, we propose the \exex{} scheme as an effective algorithm for the solution of the cardiac EFSI problem, in particular when one is focused on the ejection or filling phases or when a modular approach (i.e. the use of separate codes for the subproblems) is needed. In the latter case, we also propose the use of the \exin{2} scheme as a competitive approach for the isovolumetric phases, i.e. when one is interested in the whole heartbeat. While we have shown that the proposed loosely coupled scheme can achieve stability and efficiency in both idealized and realistic settings, further investigations are in order to determine optimal values for the Robin coefficient $\alpha$ when realistic geometries and constitutive laws are considered.

\section*{Acknowledgements}
This project has received funding from the European Research Council (ERC) under the European Union's Horizon 2020 research and innovation programme (grant agreement No 740132, iHEART - An Integrated Heart Model for the simulation of the cardiac function, P.I. Prof. A. Quarteroni).
\begin{center}
  \raisebox{-.5\height}{\includegraphics[width=.15\textwidth]{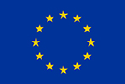}}
  \hspace{2cm}
  \raisebox{-.5\height}{\includegraphics[width=.15\textwidth]{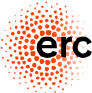}}
\end{center}

We gratefully acknowledge the CINECA award under the ISCRA initiative, for the availability of high performance computing resources and support under the project IsB25\_MathBeat, P.I. A. Quarteroni, 2021-2022.

\appendix

\section{Solid mechanics constitutive laws}
\label{app:guccione}

Given the solid displacement $\mathbf d$ and the associated deformation gradient $F = I + \grad\mathbf d$, the strain energy function associated to the Guccione constitutive law is computed as \cite{guccione1991finite,usyk2002computational}
\begin{equation*}
    \mathcal{W}_\mathrm{G}(\mathbf d) = \frac{c_\text{G}}{2}\left(\exp(Q) - 1\right) + \frac{\kappa_\text{G}}{2}(J - 1)\log(J)\,,
\end{equation*}
where
\begin{gather*}
    Q = \sum_{\mathrm{i, j} \in \{\mathrm{f}, \mathrm{s}, \mathrm{n}\}} a_\mathrm{i, j}\left(E \mathbf j\cdot \mathbf i\right)^2\;, \\
    E = \frac{1}{2}\left(F^T F - I\right)\;, \\
    J = \det F \;.
\end{gather*}
In the above equations, $c_\mathrm{G}$, $\kappa_\mathrm{G}$ and $a_\text{i,j}$ are positive parameters.

We consider for Test C (\cref{sec:zygote}) a neo-Hookean material, whose strain energy function is defined as \cite{ogden2013non}
\begin{equation*}
    \mathcal{W}_\mathrm{NH}(\mathbf d) = \frac{\mu_\text{NH}}{2}\left(J^{-\frac{2}{3}}F:F - 3\right) + \frac{\kappa_\text{NH}}{4}\left((J - 1)^2 + \log^2(J)\right)\;.
\end{equation*}
Coefficients $\mu_\mathrm{NH}$ and $\kappa_\mathrm{NH}$ are positive parameters.

\section{Model parameters}
\label{app:parameters}

\Cref{tab:parameters} reports value of physical parameters used in our model for the prolate ellipsoid test cases (Tests A and B). The parameters used for the realistic ventricle (Test C) are reported in \cref{tab:parameters-zygote}. For brevity, we only report those parameters whose values are different from the corresponding ones in Tests A and B.

\begin{table}
    \centering
    \begin{tabular}{c c S[table-format=4.3] l}
        \toprule
        & \textbf{Parameter} & \textbf{Value} & \textbf{Unit} \\
        \midrule\multirow{3}{2cm}{Electroph.}
        & $\sigma_\mathrm{m}^\mathrm{f}$ & 1.68e-4 & \si{\square\metre\per\second} \\
        & $\sigma_\mathrm{m}^\mathrm{s}$ & 0.769e-4 & \si{\square\metre\per\second} \\
        & $\sigma_\mathrm{m}^\mathrm{n}$ & 0.248e-4 & \si{\square\metre\per\second} \\
        \midrule\multirow{1}{2cm}{Force gen.}
        & $T_\text{act,max}$ & 500 & \si{\kilo\pascal} \\
        \midrule\multirow{5}{2cm}{Mechanics}
        & $\rho\solid$ & 1 & \si{\gram\per\cubic\centi\metre} \\
        & $K_\perp^\text{epi}$ & 10 & \si{\kilo\pascal\per\metre} \\
        & $K_\parallel^\text{epi}$ & 20 & \si{\kilo\pascal\per\metre} \\
        & $C_\perp^\text{epi}$ & 20 & \si{\kilo\pascal\second\per\metre} \\
        & $C_\parallel^\text{epi}$ & 2 & \si{\kilo\pascal\second\per\metre} \\
        \midrule\multirow{8}{2cm}{Guccione}
        & $c_\text{G}$ & 0.88 & \si{\kilo\pascal} \\
        & $a_\mathrm{ff}$ & 8 & \\
        & $a_\mathrm{ss}$ & 6 & \\
        & $a_\mathrm{nn}$ & 3 & \\
        & $a_\mathrm{fs}$ & 12 & \\
        & $a_\mathrm{fn}$ & 3 & \\
        & $a_\mathrm{sn}$ & 3 & \\
        & $\kappa_\text{G}$ & 50 & \si{\kilo\pascal} \\
        \midrule\multirow{5}{2cm}{Fluid}
        & $\rho\fluid$ & 1.06 & \si{\gram\per\cubic\centi\metre} \\
        & $\mu$ & 3.5e-3 & \si{\pascal\second} \\
        & $p_\text{MV}$ & 1333 & \si{\pascal} \\
        & $p_\text{AV}^0$ & 9000 & \si{\pascal} \\
        & $R_\text{AV}$ & 1.3e7 & \si{\kilo\gram\per\second\per\metre\tothe{4}} \\
        \bottomrule
    \end{tabular}
    \caption{Model parameters used in the prolate ellipsoid test cases (Tests A, B and C).}
    \label{tab:parameters}
\end{table}

\begin{table}
    \centering
    \begin{tabular}{c c S[table-format=4.3] l}
        \toprule
        & \textbf{Parameter} & \textbf{Value} & \textbf{Unit} \\
        \midrule\multirow{3}{2cm}{Electroph.}
        & $\sigma_\mathrm{m}^\mathrm{f}$ & 2e-4 & \si{\square\metre\per\second} \\
        & $\sigma_\mathrm{m}^\mathrm{s}$ & 1.05e-4 & \si{\square\metre\per\second} \\
        & $\sigma_\mathrm{m}^\mathrm{n}$ & 0.55e-4 & \si{\square\metre\per\second} \\
        \midrule\multirow{5}{2cm}{Mechanics}
        & $K_\perp^\text{epi}$ & 200 & \si{\kilo\pascal\per\metre} \\
        & $K_\parallel^\text{epi}$ & 20 & \si{\kilo\pascal\per\metre} \\
        & $C_\perp^\text{epi}$ & 20 & \si{\kilo\pascal\second\per\metre} \\
        & $C_\parallel^\text{epi}$ & 2 & \si{\kilo\pascal\second\per\metre} \\
        \midrule\multirow{2}{2cm}{Neo-Hooke}
        & $\mu_\text{NH}$ & 5000 & \si{\kilo\pascal} \\
        & $\kappa_\text{NH}$ & 5000 & \si{\kilo\pascal} \\
        \midrule\multirow{1}{2cm}{Fluid}
        & $R_\text{AV}$ & 1e7 & \si{\kilo\gram\per\second\per\metre\tothe{4}} \\
        \bottomrule
    \end{tabular}
    \caption{Test D. Model parameters used in the realistic ventricle test case. For brevity, we only report parameters whose values are different from the corresponding ones in Tests A, B and C. The latter can be found in \cref{tab:parameters}.}
    \label{tab:parameters-zygote}
\end{table}

\clearpage

\bibliographystyle{abbrv}
\bibliography{bibliography}

\end{document}